\documentclass{article}%
\usepackage{graphicx}
\usepackage{amsmath}%
\usepackage{amsfonts}%
\usepackage{amssymb}
\newtheorem{theorem}{Theorem}[section]

\newtheorem{corollary}[theorem]{Corollary}

\newtheorem{definition}[theorem]{Definition}
\newtheorem{example}[theorem]{Example}

\newtheorem{lemma}[theorem]{Lemma}

\newtheorem{proposition}[theorem]{Proposition}
\newtheorem{remark}[theorem]{Remark}

\newtheorem{theorem and definition}[theorem]{Theorem and Definition}
\newtheorem{proposition and definition}[theorem]{Proposition and Definition}
\newenvironment{proof}[1][Proof]{\textbf{#1.} }{\ \rule{0.5em}{0.5em}}

\begin{document}

\title{The Curvature and Index of Completely Positive Maps}
\author{Paul S. Muhly\thanks{Supported by grants from the U.S. National Science
Foundation and from the U.S.-Israel Binational Science Foundation.}\\Department of Mathematics\\University of Iowa\\Iowa City, IA 52242\\\texttt{muhly@math.uiowa.edu}
\and Baruch Solel\thanks{Supported by the U.S.-Israel Binational Science Foundation
and by the Fund for the Promotion of Research at the Technion}\\Department of Mathematics\\Technion\\32000 Haifa\\Israel\\\texttt{mabaruch@techunix.technion.ac.il}}
\maketitle

\begin{abstract}
We study conjugacy invariants for completely positive maps that are inspired
by the concept of curvature introduced for commuting $d$-tuples of
contractions by Arveson.

\textbf{2000 Subject Classificiation }Primary: 46L53, 46L55, 46L57, 46L87,
47L55. Secondary: 46L07, 46L08.

\end{abstract}

\section{\bigskip Introduction}

In \cite{wA00}, Arveson defined a notion of curvature for a commutative
$d$-tuple of operators $T=(T_{1},T_{2},\ldots T_{d})$, where the $T_{i}$ act
on a Hilbert space $H$ and where it is assumed that $\sum_{i=1}^{d}T_{i}%
T_{i}^{\ast}\leq I_{H}$. He showed that his concept is fundamentally an
artifact of the contractive, normal, completely positive map $\Theta_{T}$
defined on $B(H)$ by the formula%
\begin{equation}
\Theta_{T}(X)=\sum_{i=1}^{d}T_{i}XT_{i}^{\ast}\text{,} \label{cpmap}%
\end{equation}
$X\in B(H)$. Subsequently, Kribs \cite{dK01} and Popescu \cite{gP01} defined
and studied a notion of curvature for arbitrary,
\emph{not-necessarily-commuting}, $d$-tuples of contractions $T=(T_{1}%
,T_{2},\ldots T_{d})$ such that $\sum_{i=1}^{d}T_{i}T_{i}^{\ast}\leq I_{H}$.
Kribs's definition \cite[Definition 2.4]{dK01} of this curvature is directly
in terms of the map $\Theta_{T}$:
\begin{equation}
K(T):=(d-1)\lim_{k\rightarrow\infty}\frac{tr(I-\Theta_{T}^{k}(I))}{d^{k}%
}\text{.} \label{KribsDef}%
\end{equation}
(It of course needs to be proved that this limit always exists. It does, as
both Kribs and Popescu show.) Popescu's definition is different \cite[Equation
(0.1)]{gP01}; it is based on his notion of Poisson transform (see \cite{gP99})
and is closely aligned with Arveson's definition. However, he shows in Theorem
2.3 of \cite{gP01} that his definition is the same as Kribs's.\footnote{We
note in passing that $K(T)$ does not agree with Arveson's definition of
curvature when $T$ is a commutative $d$-tuple. However, Popescu shows that
Arveson's curvature can be gotten from his formula for defining $K(T)$ by
``compressing to symmetric Fock space'' \cite[Corollary 2.4]{gP01}.} Kribs and
Popescu show that the curvature $K(T)$ is a unitary invariant for $T$ that
measures its ``departure from being free'', whereby a ``free $d$-tuple'' we
mean one, $T=(T_{1},T_{2},\ldots T_{d})$, such that the $T_{i}$ are pure
isometries, i.e. multiples of the unilateral shift, with orthogonal ranges.
For such a $d$-tuple, $\Theta_{T}$ is an endomorphism of $B(H)$ and it is not
difficult to see that $K(T)$ coincides with the rank of $I-\Theta_{T}(I)$.
Further, owing to the Bunce-Frazho-Popescu dilation theorem (see \cite{jB84,
aF82, aF84, gP89}), every $d$-tuple $T$ such that $\Theta_{T}^{k}%
(I)\rightarrow0$ in the strong operator topology can be ``dilated'' to a
``free $d$-tuple'' $S$ and $K(T)\leq K(S)$, with equality holding if and only
if $T$ is itself free. In this case, $T=S$.

These results and others in \cite{dK01} and \cite{gP01}, together with
Stinespring's famous dilation theorem \cite{stinespring}, suggest the
intriguing possibility of defining the curvature of an arbitrary (contractive,
normal) completely positive map $\Theta$ on $B(H)$ via a formula like
(\ref{KribsDef}). For after all, thanks to Stinespring's analysis, every such
map is a $\Theta_{T}$ for some $d$-tuple $T$. One might hope for a rich
interplay between the new ``geometric'' invariants of $\Theta$, such as
curvature, and conjugacy or dynamical invariants of $\Theta$ as a map on
$B(H)$. \ The problem, however, is this: A contractive, normal completely
positive map on $B(H)$ can be represented by many different $d$-tuples. Even
the number $d$ is not uniquely determined. After all, the $d$-tuple
$T=(T_{1},T_{2},\ldots T_{d})$ and the $d+1$-tuple $\tilde{T}=(T_{1}%
,T_{2},\ldots T_{d},0)$ determine the same completely positive map. However,
if $K(T)$ is finite, then $K(\tilde{T})$ will be zero. So unless there is a
canonical way to represent a completely positive map on $B(H)$ in terms of a
$d$-tuple, there does not seem to be much hope in developing a notion of
curvature for completely positive maps - at least not one that proceeds along
the lines of the formula (\ref{KribsDef}).

Fortunately, however, there \emph{is }a canonical way to represent a normal,
contractive, completely positive map in terms of a $d$-tuple. This was
observed by Arveson in \cite{arvindex} and a generalization of his analysis
was developed by us to study completely positive maps on general von Neumann
algebras \cite{QMP}. It is the starting point of the present paper. Our
objective is to show how to define a concept of curvature that generalizes
Kribs's definition (\ref{KribsDef}) for \emph{any }completely positive map on
\emph{any} semifinite factor and which leads, more or less, to the same type
of results that he and Popescu found. Here, roughly and incompletely, is a
synopsis of what we do. (Full definitions and details will be given in the
body of the paper.)

Let $N$ be a semifinite factor, with a faithful normal trace $tr$, acting on a
Hilbert space $H$ and let $\Theta$ be a contractive, normal, completely
positive map on $N$. Form the Stinespring dilation $\pi:N\rightarrow
B(N\otimes_{\Theta}H)$ and let $\mathcal{E}_{\Theta}:=\{X:H\rightarrow
N\otimes_{\Theta}H\mid Xa=\pi(a)X$, $a\in N\}$. Then $\mathcal{E}_{\Theta}$
has the structure of $W^{\ast}$-correspondence over \emph{the commutant} of
$N$, $N^{\prime}$. We call $\mathcal{E}_{\Theta}$ the Arveson-Stinespring
correspondence associated with $N$. (See Proposition and Definition
\ref{esubtheta}.) This means, roughly, that $\mathcal{E}_{\Theta}$ is a
bimodule over $N^{\prime}$ and that there is $N^{\prime}$-valued inner product
on $\mathcal{E}_{\Theta}$ making $\mathcal{E}_{\Theta}$ a right Hilbert
$C^{\ast}$-module over $N^{\prime}$. The ``identity representation'' of
$\mathcal{E}$ on $H$ (Definition \ref{idenrep}) is a pair $(T,\sigma)$, where
$\sigma$ is the identity representation of $N^{\prime}$ on $H$ and where
$T:\mathcal{E}_{\Theta}\rightarrow B(H)$ is a completely contractive bimodule
map that is constructed explicitly in terms of the ingredients of the
Stinespring dilation. (The ``bimodule'' condition means that $T(a\xi
b)=aT(\xi)b$, $a,b\in N^{\prime}$ and $\xi\in\mathcal{E}_{\Theta}$.) The map
$T$, in turn, defines a Hilbert space contraction operator $\tilde
{T}:\mathcal{E}_{\Theta}\otimes_{\sigma}H\rightarrow H$ and we showed in
\cite[Corollary 2.23]{QMP} that $\Theta$ is given by the formula%
\begin{equation}
\Theta(a)=\tilde{T}(I_{\mathcal{E}_{\Theta}}\otimes a)\tilde{T}^{\ast}\text{,}
\label{tt}%
\end{equation}
$a\in N$ (see Proposition \ref{equality} below, also).

Now the correspondence $\mathcal{E}_{\Theta}$ has a natural dimension $d$,
which is a non-negative real number or $+\infty$ (Definition \ref{dime}). This
dimension is defined in terms of the trace $tr$ on $N$, but it is independent
of how $N$ is represented, so long as the commutant of the representation is
finite (Theorem \ref{uniqued}). Because of its invariance under
representations of $N$, we call the dimension $d$ the index of $\Theta$ and
write $d=d(\Theta)$ (Definition \ref{index}). If $\Theta$ happens to be given
by an $n$-tuple of operators $(t_{1},t_{2},\ldots t_{n})$ in $N$, i.e., if
$\Theta(a)=\sum_{i=1}^{n}t_{i}at_{i}^{\ast}$, $a\in N$, then $d(\Theta)$ turns
out to be the vector space dimension of the complex linear span of the $t_{i}$
in $N$ (Proposition \ref{ti}). In that case, too, the identity representation
$(T,\sigma)$ of $\mathcal{E}_{\Theta}$ gives a canonical $d$-tuple $(\tilde
{t}_{1},\tilde{t}_{2},\ldots\tilde{t}_{d})$ of operators in $N$ representing
$\Theta$ through equation (\ref{tt}), i.e., such that $\Theta(a)=\sum
_{i=1}^{d}\tilde{t}_{i}a\tilde{t}_{i}^{\ast}$, $a\in N$.

If $d=d(\Theta)$ is finite, then the limit%
\[
\lim_{k\rightarrow\infty}\frac{tr(I-\Theta^{k}(I))}{\sum_{j=0}^{k-1}d^{j}}%
\]
exists as a positive real number or $+\infty$. This limit is our definition of
the curvature $K(\Theta,tr)$ of $\Theta$ (Definition \ref{curvcp}). Evidently,
$K(\Theta,tr)$ depends on the trace $tr$. However, the quantity $K(\Theta
,tr)/tr(I-\Theta(I))$ doesn't and truly should be thought of as the curvature
of $\Theta$; we call this quantity the normalized curvature of $\Theta$
(Definition \ref{k1}) and show that it is a conjugacy invariant of $\Theta$
(Theorem \ref{invariant}). When the index of $\Theta$, $d$, is strictly larger
than $1$, then the limit defining $K(\Theta,tr)$ turns out to be the same as
$(d-1)\lim_{k\rightarrow\infty}\frac{tr(I-\Theta^{k}(I))}{d^{k}}$, which is
Kribs's definition (\ref{KribsDef}) in the case when $N=B(H)$.

In \cite{QMP}, we showed how to dilate the completely positive map $\Theta$ on
$N$ to an \emph{endomorphism} $\alpha$ of a von Neumann algebra $R$ in which
$N$ sits as a corner. This was based on our dilation theory for
representations of correspondences, which, in turn, generalizes the
Bunce-Frazho-Popescu dilation theorem (see \cite[Theorem 3.3]{MStensor} and
\cite[Theorem and Definition 2.18]{QMP}). We show in Theorem \ref{invofdil}
that if $d(\Theta)$ is finite, then so is $d(\alpha)$ and the two are equal.
Further, $K(\alpha,tr_{R^{^{\prime}}})\geq K(\Theta,tr_{N^{^{\prime}}})$. The
curvature of the endomorphism $\alpha$ is $tr_{R^{\prime}}(I-\alpha(I))$. If
$\Theta$ is pure, in the sense that the sequence $\{\Theta^{k}(I)\}_{k\geq0}$
converges strongly to zero, then so is $\alpha$ and in this case
$K(\alpha,tr_{R^{^{\prime}}})=K(\Theta,tr_{N^{^{\prime}}})$ if and only if
$\Theta=\alpha$. Thus, given that endomorphisms may be expressed in terms of
isometric representations of correspondences and that there is an analogue of
``pure isometric representation'' (see \cite{wold} and Definition \ref{pure}
below), we find that our notion of curvature for a completely positive map may
also be viewed a measure of how much the map deviates from a ``free''
representation. \ Thus, our notion of curvature for contractive, normal,
completely positive maps on semifinite factors captures the salient features
of the concept defined in equation (\ref{KribsDef}).

The constructs we have defined are easy to calculate in some situations. We
offer a number of examples in Section 3. Of particular note is Example
\ref{singlet} in which we show that if $N$ is a semifinite factor with
normalized trace $tr$ and if $\Theta(a)=tat^{\ast}$ for a contraction $t\in N$
such that $tr(I-tt^{\ast})<\infty$, then $d(\Theta)=1$ and $K(\Theta
,tr)=tr(I-tt^{\ast})-tr((I-t^{\ast}t)^{1/2}(I-t_{\infty})(I-t^{\ast}t)^{1/2}%
)$. This generalizes work of Parrott and Levy \cite{Parrott, rLpp00}.

The next section is devoted to assembling material which we will use, from
various sources. We discuss the general theory of $W^{\ast}$-correspondences,
their representation and dilation theory, and the notion of dimension for
correspondences over semifinite factors. In the third section, we develop the
notion of the curvature of a representation of a correspondence. We use this
material in the last section to define the curvature of a completely positive
map through the curvature of the identity representation of its Stinespring
correspondence. We show that the constructs we define are independent of any
representation and so are intrinsic features of completely positive maps. \ We
also calculate a number of examples.

\section{Preliminaries\label{prelim}}

In this section we collect or develop a variety of facts that will be used in
the sequel. We organize them into several subsections which are somewhat
disjoint as presented here, but which will be blended in the next two sections.

\subsection{Dimensions of Representations of Semifinite Factors}

Let $M$ be a semifinite factor with normal, faithful semifinite trace $\tau$.
We do not preclude the possibility that $M$ is finite. However, when $M$ is
finite, we do not necessarily assume that the trace $\tau$ is normalized so
that $\tau(I)=1$. \ We will make explicit our assumptions on the normalization
of traces as they arise. We write $L^{2}(M)$, or $L^{2}(M,\tau)$, for the
space obtained by the GNS construction applied to $M$ and $\tau$. Every
element $a$ of $M$ then defines an operator of left multiplication on
$L^{2}(M)$ , denoted $\lambda(a)$ , and an operator of right multiplication,
denoted $\rho(a)$. The maps $\lambda$ and $\rho$ are $\ast$-representations of
$M$ and the opposite algebra of $M$, $M^{opp}$, respectively. The symmetry
between the left and right representations of $M$ is implemented by a
conjugate linear isometry $J:L^{2}(M)\rightarrow L^{2}(M)$, which is simply
the extension to $L^{2}(M)$ of the map $x\rightarrow x^{\ast}$ defined on $M$.
We have $J\lambda(M)J=\rho(M)$.

Suppose $H$ is a (left) $M$-module; i.e. suppose that there is a unital,
normal, $\ast$-representation $\sigma$ of $M$ on $H$.\footnote{In this paper,
we consider only normal $\ast$-representations of von Neumann algebras.
\ Also, except in certain cases involving not-necessarily unital
endomorphisms, which will be clearly identified, all our representations will
be unital. \ Consequently, we often drop the adjectives ``normal'', ``$\ast
$-'', and ``unital'' and speak simply of ``representations''.} It is known
that there is an $M$-linear isometry $u:H\rightarrow L^{2}(M)\otimes K$, where
$K$ is a separable, infinite dimensional Hilbert space ; i.e. $u\sigma
(a)=(\lambda(a)\otimes I_{K})u$ for every $a$ in $M$. The projection
$p:=uu^{\ast}$ is in the commutant of $\lambda(M)\otimes I_{K}$, and we find
that $u\sigma(M)^{\prime}u^{\ast}=p(\lambda(M)\otimes I_{K})^{\prime
}p\subseteq(\lambda(M)\otimes I_{K})^{\prime}$. Note that the commutant of
$\lambda(M)\otimes I_{K}$ is the semifinite factor $\rho(M)\otimes B(K)$ so
that every element of this algebra can be written as a matrix $(\rho(a_{ij}%
))$. Given $x$ in $\sigma(M)_{+}^{\prime}$ write $uxu^{\ast}$ as $(\rho
(a_{ij}))$ and define $tr_{\sigma(M)^{\prime}}(x)$ (or, sometimes,
$tr_{\sigma(M)^{\prime},H}(x)$) to be $\sum\tau(a_{ii})$. This yields a
normal, semifinite, faithful trace, also denoted $tr_{\sigma(M)^{\prime}}$ or
$tr_{\sigma(M)^{\prime},H}$ on the algebra $\sigma(M)^{\prime}$.

\begin{definition}
\label{dimH}Let $M$ be a semifinite factor with a prescribed normal semifinite
trace $tr$.

\begin{itemize}
\item[(1)] Given a representation $\sigma$ of $M$ on a Hilbert space $H$, the
\emph{natural trace} on $\sigma(M)^{\prime}$ is the trace $tr_{\sigma
(M)^{\prime}}$ just defined.

\item[(2)] The (left) dimension of $H$ (as an $M$-module) is
\[
dim_{M}(H)=tr_{\sigma(M)^{\prime}}(p)
\]
where $p$ is as above.
\end{itemize}
\end{definition}

\begin{lemma}
\label{Ne}Let $M$ be a semifinite factor represented (necessarily faithfully)
on a Hilbert space $H$, with commutant $M^{\prime}$. Let $e$ be a projection
in $M^{\prime}$ and let $x$ be a positive element of $eM^{\prime}e$. Then
\begin{equation}
tr_{M^{\prime}}(x)=tr_{(eM)^{\prime},eH}(x)
\end{equation}
\end{lemma}

\begin{proof}
As discussed above, we have an $M$-linear isometry $u$ from $H$ to
$L^{2}(M)\otimes K$ whose range is a projection $p$ in $(\lambda(M)\otimes
I_{K})^{\prime}$. The positive element $uxu^{\ast}$ also lies in
$(\lambda(M)\otimes I_{K})^{\prime}$ and can be written as a matrix
$(\rho(x_{ij}))$ for $x_{ij}$ in $M$. Note that the map $a\rightarrow ae$ is a
representation of $M$ on $eH$. We write $Q$ for the image of this
representation, $Me$. If $\tau$ is the trace on $M$, then we transport it to
one, $\tau_{0}$, on $Q$ via the formula $\tau_{0}(ae)=\tau(a)$. The map
$V:L^{2}(M)\rightarrow L^{2}(Q)$ extending the map $a\mapsto ae$ is a Hilbert
space isomorphism and $\lambda_{Q}(ae)=V\lambda_{M}(a)V^{\ast}$. Then the
composition of $ue$ with the map $V\otimes I_{K}$ is a $Q$-module isometry
from $eH$ to $L^{2}(Q)\otimes K$ and we denote it by $v$. A straightforward
computation shows that $vxv^{\ast}=(V\otimes I)uxu^{\ast}(V^{\ast}\otimes I)$
has the matrix representation $(\rho(x_{ij}e))$ and, therefore,
$tr_{(eM)^{\prime},eH}(x)=\sum\tau_{0}(x_{ij}e)=\sum\tau(x_{ij})=tr_{M^{\prime
}}(x)$.
\end{proof}

\subsection{$W^{\ast}$-correspondences and their dimensions}

We begin by recalling the notion of a $W^{\ast}$-correspondence. For the
general theory of Hilbert $C^{\ast}$-modules which we use, we will follow
\cite{lance}. In particular, a Hilbert $C^{\ast}$-module will be a
\emph{right} Hilbert $C^{\ast}$-module.

\begin{definition}
\label{correspondence}Let $M$ and $N$ be von Neumann algebras and let
$\mathcal{E}$ be a (right) Hilbert $C^{\ast}$-module over $N$. Then
$\mathcal{E}$ is called a \emph{Hilbert $W^{\ast}$-module} over $N$ in case it
is self dual (i.e. every continuous $N$-module map from $\mathcal{E}$ to $N$
is implemented by an element of $\mathcal{E}$). It is called a \emph{$W^{\ast
}$-correspondence} from $M$ to $N$ if it is also endowed with a structure of a
left $M$-module via a normal $\ast$-homomorphism $\varphi:M\rightarrow
\mathcal{L}(\mathcal{E})$.(Here $\mathcal{L}(\mathcal{E})$ is the algebra of
all bounded, adjointable, module maps on $\mathcal{E}$. For a Hilbert
$W^{\ast}$-module it is known to be a von Neumann algebra). A \emph{\ $W^{\ast
}$-correspondence over $M$} is simply a $W^{\ast}$-correspondence from $M$ to
$M$.
\end{definition}

If $\mathcal{E}$ is a $W^{\ast}$-correspondence from $M$ to $N$ and if
$\mathcal{F}$ is a $W^{\ast}$-correspondence from $N$ to $Q$, then the
balanced tensor product, $\mathcal{E}\otimes_{N}\mathcal{F}$ is a $W^{\ast}%
$-correspondence from $M$ to $Q$. It is defined as the Hausdorff completion of
the algebraic balanced tensor product with the internal inner product given
by
\[
\langle\xi_{1}\otimes\eta_{1},\xi_{2}\otimes\eta_{2}\rangle=\langle\eta
_{1},\varphi(\langle\xi_{1},\xi_{2}\rangle_{\mathcal{E}})\eta_{2}%
\rangle_{\mathcal{F}}%
\]
for all $\xi_{1}$ , $\xi_{2}$ in $\mathcal{E}$ and $\eta_{1}$ , $\eta_{2}$ in
$\mathcal{F}$. The left and right actions are defined by
\[
\varphi_{\mathcal{E}\otimes_{N}\mathcal{F}}(a)(\xi\otimes\eta)b=\varphi
_{\mathcal{E}}(a)\xi\otimes\eta b
\]
for $a$ in $M$, $b$ in $Q$, $\xi$ in $\mathcal{E}$ and $\eta$ in $\mathcal{F}$.

Given a $W^{\ast}$-correspondence $\mathcal{E}$ over $M$, the \emph{full Fock
space} over $\mathcal{E}$ will be denoted by $\mathcal{F}(\mathcal{E})$, so
$\mathcal{F}(\mathcal{E})=M\oplus\mathcal{E}\oplus\mathcal{E}^{\otimes2}%
\oplus\cdots$. It is also a $W^{\ast}$-correspondence over $M$ with left
action $\varphi_{\infty}$ (or $\varphi_{\mathcal{E},\infty}$) given by the
formula
\[
\varphi_{\infty}(a)=diag(a,\varphi(a),\varphi^{(2)}(a),\cdots)\text{,}%
\]
where $\varphi^{(n)}(a)(\xi_{1}\otimes\xi_{2}\otimes\cdots\xi_{n}%
)=(\varphi(a)\xi_{1})\otimes\xi_{2}\otimes\cdots\xi_{n}$ . For $\;\xi
\in\mathcal{E}\;$ we write $T_{\xi}$ for the creation operator on
$\mathcal{F}(\mathcal{E})$ : $T_{\xi}\eta=\xi\otimes\eta,\;\eta\in
\mathcal{F}(\mathcal{E})\;$.

We also recall that any $W^{\ast}$-correspondence over $M$ carries a natural
weak topology, called the $\sigma$-topology (see \cite{sigmatop}).This is the
topology defined by the functionals $\;f(\cdot)=\sum_{n=1}^{\infty}\omega
_{n}(\langle\eta_{n},\cdot\rangle)\;$ where the $\eta_{n}$ lie in
$\mathcal{E}$ , the $\omega_{n}$ lie in $M_{\ast}$ , and $\sum\Vert\omega
_{n}\Vert\Vert\eta_{n}\Vert<\infty$.

We shall need some of the concepts and results of Jones' index theory and we
will refer to \cite{J} or \cite{JS} for the basic results.

\begin{definition}
\label{dim} If $M$ is a semifinite factor, a Hilbert space $H$ is said to be
an \emph{$M$-$M$ bimodule} if $H$ is a left $M$-module (whose structure is
given by a unital normal representation $\pi_{l}$ of $M$ on $H$), $H$ is a
right $M$-module (whose structure is given by a unital normal representation
$\pi_{r}^{0}$ of $M^{opp}$ or a unital normal antirepresentation $\pi_{r}$ of
$M$ on $H$) and the actions commute (i.e. $\pi_{r}(M)\subseteq(\pi
_{l}(M))^{\prime}$). The dimensions of $H$ with respect to $\pi_{l}$ and
$\pi_{r}$ will be denoted $dim_{M-}(H)$ and $dim_{-M}(H)$ respectively. We
shall call the $M$-bimodule \emph{bifinite} if both of these numbers are finite.
\end{definition}

Observe that if $H$ is a left $M$-module, it can be viewed as a $W^{\ast}%
$-correspondence from $M$ to $\mathbb{C}$. If $\mathcal{E}$ is a $W^{\ast}%
$-correspondence over $M$ then the balanced tensor product $\mathcal{E}%
\otimes_{M}H$ is also a left $M$-module. In particular, if $H$ is the $M-M$
bimodule $L^{2}(M)$ then the tensor product is an $M-M$ bimodule. \ The map
$\mathcal{E}\rightarrow\mathcal{E}\otimes_{M}L^{2}(M)$ defines a bijection
between $W^{\ast}$-correspondences over $M$ and $M-M$ bimodules. This map is
explored in \cite{sigmatop}.

\begin{definition}
\label{dime}If $\mathcal{E}$ is a $W^{\ast}$-correspondence over a semifinite
factor $M$, then the \emph{left dimension } of $\mathcal{E}$ is defined to be
$dim_{M-}(\mathcal{E}\otimes L^{2}(M))$ and will be written $dim_{l}%
(\mathcal{E})$. Similarly, the \emph{\ right dimension} of $\mathcal{E}$ is
$dim_{-M}(\mathcal{E}\otimes L^{2}(M))$ and will be written $dim_{r}%
(\mathcal{E})$. A $W^{\ast}$-correspondence is said to be \emph{left-finite}
(respectively, \emph{right-finite}) if $dim_{l}(\mathcal{E})$ (respectively,
$dim_{r}(\mathcal{E})$) is finite. It is said to be \emph{bifinite} if both
these numbers are finite.
\end{definition}

Let $\mathcal{E}$ be a $W^{\ast}$-correspondence over the von Neumann algebra
$M$ and let $H$ be a left $M$-module, with associated normal representation
$\sigma$. Then there is an \emph{induced representation} $\sigma^{\mathcal{E}%
}:\mathcal{L}(\mathcal{E})\rightarrow B(\mathcal{E}\otimes_{\sigma}H)$ defined
by the formula $\sigma^{\mathcal{E}}(S)=S\otimes I$ \cite[Lemma 3.4]%
{MStensor}. It is not hard to see that $\sigma^{\mathcal{E}}$ is a normal
representation. By \cite[Theorem 6.23]{rieffel}, the image of the induced
representation is the commutant of the algebra of all operators of the form
$I_{\mathcal{E}}\otimes T$, where $T$ lies in the commutant of $\sigma(M)$. We
have
\[
\mathcal{L}(\mathcal{E})\otimes I_{H}=(I_{\mathcal{E}}\otimes\sigma
(M)^{\prime})^{\prime}=(I_{\mathcal{E}}\otimes\kappa(M))^{\prime},
\]
where $\kappa$ is the representation of $\sigma(M)^{\prime}$ given by the
formula $T\rightarrow I_{\mathcal{E}}\otimes T$, $T\in\sigma(M)^{\prime}$. In
the special case when $H=L^{2}(M)$, so that $\sigma=\lambda$ and $\kappa=\rho
$, we write $\pi_{l}$ and $\pi_{r}$ for the representation and the
antirepresentation that define the left and right actions of $M$ on
$\mathcal{E}\otimes L^{2}(M)$. We shall also, on occasion, write $\pi_{l}$ as
$\varphi(\cdot)\otimes I$ and $\pi_{r}$ as $I_{\mathcal{E}}\otimes\rho(\cdot)$.

\begin{lemma}
\label{L(E)}Let $\mathcal{E}$ be a bifinite $W^{\ast}$-correspondence over the
type $II_{1}$ factor $M$. Then

\begin{enumerate}
\item[(1)] The von Neumann algebra $\mathcal{L}(\mathcal{E})$ is a type
$II_{1}$ factor and $\varphi(M)$ is a subfactor.

\item[(2)] If $\sigma$ is a representation of $M$ on $H$ then
\[
\lbrack\mathcal{L}(\mathcal{E})\otimes I_{H}:\varphi(M)\otimes I_{H}%
](=[\sigma^{\mathcal{E}}(\mathcal{L}(\mathcal{E})):\sigma^{\mathcal{E}%
}(\varphi(M))])=[\mathcal{L}(\mathcal{E}):\varphi(M)].
\]

\item[(3)] $[\mathcal{L}(\mathcal{E}):\varphi(M)]=dim_{l}(\mathcal{E}%
)dim_{r}(\mathcal{E})<\infty$
\end{enumerate}
\end{lemma}

\begin{proof}
Suppose $\sigma$ is a normal representation of $M$ on $H$, as in part (2).
Since $M$ is a factor, $\sigma$ is an isomorphism and it is easy to check that
then $\sigma^{\mathcal{E}}$ is an isomorphism, so that $\mathcal{L}%
(\mathcal{E})$ is isomorphic to $\mathcal{L}(\mathcal{E})\otimes I_{H}$. In
particular, for $H=L^{2}(M)$, $\mathcal{L}(\mathcal{E})$ is isomorphic to
$\mathcal{L}(\mathcal{E})\otimes I_{L^{2}(M)}$ and that algebra is the
commutant of $I_{\mathcal{E}}\otimes\rho(M)$. As $\rho(M)$ is a type $II_{1}$
factor and $I_{\mathcal{E}}\otimes\rho(M)$ is isomorphic to it, $\mathcal{L}%
(\mathcal{E})$ is a type $II$ factor. Since $I_{\mathcal{E}}\otimes\rho(M)$ is
a finite factor and $dim_{r}(\mathcal{E})$ is finite it follows that
$\mathcal{L}(\mathcal{E})\otimes I_{L^{2}(M)}$ (and, hence also $\mathcal{L}%
(\mathcal{E})$) is finite (\cite[Proposition 2.2.6(iii)]{JS}). This proves
(1). Part (2) follows from the fact that $\sigma^{\mathcal{E}}$ is an
isomorphism and part (3) is in \cite[Discussion following Corollary 2.3.6]{JS}.
\end{proof}

\begin{lemma}
\label{trace}Let $M$ be a finite factor and let $\sigma$ be a normal
representation of $M$ on $H$. If $\mathcal{E}$ is a left-finite $W^{\ast}%
$-correspondence over $M$, then for every positive element $x$ in
$\sigma(M)^{\prime}$, we have
\[
tr_{\sigma(M)^{\prime}}(x)\cdot dim_{l}(\mathcal{E})=tr_{(\varphi(M)\otimes
I_{H})^{\prime}}(I_{\mathcal{E}}\otimes x)\text{.}%
\]
\end{lemma}

\begin{proof}
If $M=M_{n}(\mathbb{C})$, it is easy to check that there is a projection $e$
in $\sigma(M)^{\prime}$ with $tr_{\sigma(M)^{\prime}}(e)=1/n$. If $M$ is of
type $II_{1}$, the values of the trace $tr_{\sigma(M)^{\prime}}$ on
projections of $\sigma(M)^{\prime}$ form an interval containing $0$. In any
case, one can always find a projection $e$ in $\sigma(M)^{\prime}$ and a
positive integer $m$ such that $tr_{\sigma(M)^{\prime}}(e)=1/m$. Fix such a
projection. Then $dim_{\sigma(M)}eH=1/m\;$ and
\[
\sum_{i=1}^{m}\oplus eH\cong L^{2}(M)
\]
as left $M$-modules (where the action of $M$ on $eH$ is by $\sigma$ and on
$L^{2}(M)$ by $\lambda$). Tensoring by $\mathcal{E}$, we get
\[
\sum_{i=1}^{m}\oplus(\mathcal{E}\otimes eH)\cong\mathcal{E}\otimes L^{2}(M).
\]
It follows that
\[
tr_{(\rho(M)\otimes I)^{\prime}}(I_{\mathcal{E}}\otimes e)=dim_{\rho(M)\otimes
I}\mathcal{E}\otimes eH=\frac{1}{m}dim_{l}(\mathcal{E})=dim_{l}(\mathcal{E}%
)\cdot tr_{\sigma(M)^{\prime}}(e).
\]
Writing $\tau_{1}(x)=tr_{(\rho(M)\otimes I)^{\prime}}(I_{\mathcal{E}}\otimes
x)\;$ for $0\leq x\in\sigma(M)^{\prime}$,$\;$ we get a faithful, normal,
semifinite trace on $\sigma(M)^{\prime}$. Thus it is a multiple of
$tr_{\sigma(M)^{\prime}}$. The computation above shows that the multiple is
$dim_{l}(\mathcal{E})$.
\end{proof}

\begin{corollary}
\label{multdim}Let $M$ be a finite factor and let $\mathcal{E}$ and
$\mathcal{F}$ be two left-finite $W^{\ast}$-correspondences over $M$. Then
\[
dim_{l}(\mathcal{E}\otimes_{M}\mathcal{F})=dim_{l}(\mathcal{E})dim_{l}%
(\mathcal{F})\text{.}%
\]
In particular, $dim_{l}(\mathcal{E}^{\otimes n})=(dim_{l}(\mathcal{E}))^{n}$
for all positive integers $n$.
\end{corollary}

\begin{proof}
The corollary is proved by applying Lemma~\ref{trace} several times. We write
$H$ for the space $\mathcal{F}\otimes L^{2}(M)$ and $\sigma$ will be the
representation $\varphi_{\mathcal{F}}(\cdot)\otimes I_{L^{2}(M)}$ of $M$ on
$H$ (we write $\varphi_{\mathcal{E}}$, $\varphi_{\mathcal{F}}$ and
$\varphi_{\mathcal{E}\otimes\mathcal{F}}$ for the left action maps on
$\mathcal{E}$, $\mathcal{F}$ and $\mathcal{E}\otimes\mathcal{F}$
respectively). Note that $\varphi_{\mathcal{E}}(M)\otimes I_{H}=\varphi
_{\mathcal{E}\otimes\mathcal{F}}(M)\otimes I_{L^{2}(M)}$. Thus we have, for a
positive element $x$ in $\lambda(M)^{\prime}=\rho(M)$,
\[
tr_{(\varphi_{\mathcal{E}\otimes\mathcal{F}}(M)\otimes I_{L^{2}(M)})^{\prime}%
}(I_{\mathcal{E}\otimes\mathcal{F}}\otimes x)=tr_{(\varphi_{\mathcal{E}%
}(M)\otimes I_{H})^{\prime}}(I_{\mathcal{E}}\otimes I_{\mathcal{F}}\otimes x)
\]%
\[
=dim_{l}(\mathcal{E})\cdot tr_{\sigma(M)^{\prime}}(I_{\mathcal{F}}\otimes
x)=dim_{l}(\mathcal{E})dim_{l}(\mathcal{F})\cdot tr_{\lambda(M)^{\prime}}(x)
\]
But, using Lemma~\ref{trace} again, this will be equal to $tr_{\lambda
(M)^{\prime}}(x)dim_{l}(\mathcal{E}\otimes\mathcal{F})$. Now set $x=I$ to
complete the proof.
\end{proof}

\subsection{Representations of correspondences and completely positive maps}

In this subsection discuss representations of $W^{\ast}$-correspondences and
the completely positive maps associated with such representations. For more
details, please refer to \cite{MStensor} and \cite{QMP}.

\begin{definition}
Let $\mathcal{E}$ be a $W^{\ast}$-correspondence over a von Neumann algebra
$N$ and let $H$ be a Hilbert space.

\begin{enumerate}
\item[(1)] A \emph{completely contractive covariant representation} of
$\mathcal{E}$ (or, simply, a \emph{representation} of $\mathcal{E}$) in $B(H)$
is a pair $(T,\sigma)$, where

\begin{enumerate}
\item[(a)] $\sigma$ is a normal representation of $N$ in $B(H)$.

\item[(b)] $T$ is a linear, completely contractive map from $\mathcal{E}$ to
$B(H)$ that is continuous with respect to the $\sigma$-topology of
\cite{sigmatop} on $\mathcal{E}$ and the $\sigma$-weak topology on $B(H)$.

\item[(c)] $T$ is a bimodule map in the sense that $T(\varphi(a)\xi
b)=\sigma(a)T(\xi)\sigma(b),$ $\xi\in\mathcal{E}$, and $a,b\;\in\;N$.
\end{enumerate}

\item[(2)] A completely contractive covariant representation $(T,\sigma)$ of
$\mathcal{E}$ in $B(H)$ is called \emph{isometric} in case
\[
T(\xi)^{\ast}T(\eta)=\sigma(\langle\xi,\eta\rangle),
\]
for all $\xi,\eta$ in $\mathcal{E}$.
\end{enumerate}
\end{definition}

The theory developed in \cite{MStensor} applies here to prove that if a
representation $(T,\sigma)$ of $\mathcal{E}$ is given, then it determines a
contraction $\tilde{T}:\mathcal{E}\otimes_{\sigma}H\rightarrow H$ defined by
the formula
\[
\tilde{T}(\xi\otimes h)=T(\xi)h.
\]
Moreover, for every $a$ in $N$ we have
\begin{equation}
\tilde{T}(\varphi(a)\otimes I)(=\tilde{T}\sigma^{\mathcal{E}}(\varphi
(a)))=\sigma(a)\tilde{T}\text{,} \label{Ttilde}%
\end{equation}
i.e., $\tilde{T}$ intertwines $\sigma$ and $\sigma^{\mathcal{E}}\circ\varphi$.
In fact, it is shown in \cite{MStensor} that there is a bijection between
representations $(T,\sigma)$ of $\mathcal{E}$ and intertwining operators
$\tilde{T}$ of $\sigma$ and $\sigma^{\mathcal{E}}\circ\varphi$.

It is also shown in \cite{MStensor} that $(T,\sigma)$ is isometric if and only
if $\tilde{T}$ is an isometry.

\begin{definition}
A representation $(T,\sigma)$ of a $W^{\ast}$-correspondence $\mathcal{E}$ is
called \emph{\ fully coisometric} if $\tilde{T}$ is a coisometry (i.e. if
$\tilde{T}\tilde{T}^{\ast}=I_{H}$).
\end{definition}

Associated with every representation $(T,\sigma)$ we can define a completely
positive, contractive, map $\Theta=\Theta_{T}$ on the \emph{commutant} of
$\sigma(N)$. It is given in the following proposition. The proof and more
details can be found in \cite[Proposition 2.21]{QMP}.

\begin{proposition and definition}
\label{theta}Let $N$ be a von Neumann algebra, let $\mathcal{E}$ be a
$W^{\ast}$-correspondence over $N$ and let $(T,\sigma)$ be a completely
contractive covariant representation of $\mathcal{E}$ on a Hilbert space $H$.
For $S\in\sigma(N)^{\prime}$ , set
\[
\Theta(S)=\Theta_{T}(S)=\tilde{T}(I_{\mathcal{E}}\otimes S)\tilde{T}^{\ast}.
\]
Then $\Theta$ is a contractive, normal completely positive map from
$\sigma(N)^{\prime}$ into itself. It is unital if and only if $(T,\sigma)$ is
fully coisometric and it is multiplicative (i.e. a $\ast$-endomorphism of
$\sigma(N)^{\prime}$) if $(T,\sigma)$ is isometric.
\end{proposition and definition}

In addition to $\tilde{T}$ we also define the maps $\;\tilde{T}_{n}%
:\mathcal{E}^{\otimes n}\otimes H\rightarrow H\;$ by $\;\tilde{T}_{n}(\xi
_{1}\otimes\ldots\otimes\xi_{n}\otimes h)=T(\xi_{1})\cdots T(\xi_{n})h\;$ and
then we have $\;\tilde{T}_{n+m}=\tilde{T}_{n}(I_{n}\otimes\tilde{T}%
_{m})=\tilde{T}_{m}(I_{m}\otimes\tilde{T}_{n})$, where $I_{n}$ is the identity
map on $\mathcal{E}^{\otimes n}$ \cite{wold}. It follows that
\[
\Theta_{T}^{n}(S)=\tilde{T}_{n}(I_{n}\otimes S)\tilde{T}_{n}^{\ast}%
\]
for $S$ in $\sigma(M)^{\prime}$.

\begin{proposition}
\label{estimate}Let $\mathcal{E}$ be a left-finite $W^{\ast}$-correspondence
over a finite factor $M$ and let $(T,\sigma)$ be a representation of
$\mathcal{E}$ in $B(H)$. Then for every positive $x$ in $\sigma(M)^{\prime}$
we have
\[
tr_{\sigma(M)^{\prime}}(\Theta_{T}(x))\leq\Vert\tilde{T}\Vert^{2}%
dim_{l}(\mathcal{E})tr_{\sigma(M)^{\prime}}(x).
\]
\end{proposition}

\begin{proof}
Let $\tilde{\sigma}$ be the representation of $M$ on $H\oplus(\mathcal{E}%
\otimes_{\sigma}H)$ defined by the formula $\tilde{\sigma}(a)=\left(
\begin{array}
[c]{cc}%
\sigma(a) & 0\\
0 & \varphi(a)\otimes I_{H}%
\end{array}
\right)  $ and let $\mathcal{B}$ be the image of $\tilde{\sigma}$:
\begin{equation}
\mathcal{B}=\tilde{\sigma}(M)=\left\{  \left(
\begin{array}
[c]{cc}%
\sigma(a) & 0\\
0 & \varphi(a)\otimes I_{H}%
\end{array}
\right)  \in B(H\oplus(\mathcal{E}\otimes H)):\;\;a\in M\right\}  .
\label{diag1}%
\end{equation}
Write $tr_{\mathcal{B}^{\prime}}$ for the natural trace on $\mathcal{B}%
^{\prime}$ (Definition \ref{dimH}).

For $y$ in $\sigma(M)^{\prime}$ and $z$ in $(\varphi(M)\otimes I)^{\prime}$ it
is easy to check that $\left(
\begin{array}
[c]{cc}%
y & 0\\
0 & 0
\end{array}
\right)  $ and $\left(
\begin{array}
[c]{cc}%
0 & 0\\
0 & z
\end{array}
\right)  $ lie in $\mathcal{B}^{\prime}$ and that the equations,
\[
tr_{\mathcal{B}^{\prime}}\left(
\begin{array}
[c]{cc}%
y & 0\\
0 & 0
\end{array}
\right)  =tr_{\sigma(M)^{\prime}}(y)
\]
and%
\[
tr_{\mathcal{B}^{\prime}}\left(
\begin{array}
[c]{cc}%
0 & 0\\
0 & z
\end{array}
\right)  =tr_{(\varphi(M)\otimes I)^{\prime}}(z)\text{,}%
\]
are valid.

Note also that equation (\ref{Ttilde}) implies that $\left(
\begin{array}
[c]{cc}%
0 & \tilde{T}\\
0 & 0
\end{array}
\right)  $ lies in $\mathcal{B}$. We therefore have, for a positive element
$x=c^{\ast}c$ in $\sigma(M)^{\prime}$, that
\[
tr_{\sigma(M)^{\prime}}(\Theta_{T}(x))=tr_{\sigma(M)^{\prime}}(\tilde
{T}(I\otimes x)\tilde{T}^{\ast})=tr_{\mathcal{B}^{\prime}}\left(
\begin{array}
[c]{cc}%
\tilde{T}(I\otimes x)\tilde{T}^{\ast} & 0\\
0 & 0
\end{array}
\right)
\]%
\[
=tr_{\mathcal{B}^{\prime}}\left(
\begin{array}
[c]{cc}%
0 & \tilde{T}(I\otimes c^{\ast})\\
0 & 0
\end{array}
\right)  \left(
\begin{array}
[c]{cc}%
0 & 0\\
(I\otimes c)\tilde{T}^{\ast} & 0
\end{array}
\right)
\]%
\[
=tr_{\mathcal{B}^{\prime}}\left(
\begin{array}
[c]{cc}%
0 & 0\\
(I\otimes c)\tilde{T}^{\ast} & 0
\end{array}
\right)  \left(
\begin{array}
[c]{cc}%
0 & \tilde{T}(I\otimes c^{\ast})\\
0 & 0
\end{array}
\right)
\]%
\[
=tr_{\mathcal{B}^{\prime}}\left(
\begin{array}
[c]{cc}%
0 & 0\\
0 & (I\otimes c)\tilde{T}^{\ast}\tilde{T}(I\otimes c^{\ast})
\end{array}
\right)  =tr_{(\varphi(M)\otimes I)^{\prime}}((I\otimes c)\tilde{T}^{\ast
}\tilde{T}(I\otimes c^{\ast}))
\]%
\[
\leq\Vert\tilde{T}\Vert^{2}tr_{(\varphi(M)\otimes I)^{\prime}}(I\otimes
x)=\Vert\tilde{T}\Vert^{2}dim_{l}(\mathcal{E})tr_{\sigma(M)^{\prime}%
}(x)\text{,}%
\]
where the last equality follows from Lemma~\ref{trace}.
\end{proof}

\begin{corollary}
\label{isomest}Let $\mathcal{E}$ be a left-finite $W^{\ast}$-correspondence
over a finite factor $M$ and let $(T,\sigma)$ be a representation of
$\mathcal{E}$, as in Proposition~\ref{estimate}. Write $\;\Delta=(I-\tilde
{T}^{\ast}\tilde{T})^{1/2}\;$. Then, for every positive element $x$ in
$\sigma(M)^{\prime}$,
\[
dim_{l}(\mathcal{E})\cdot tr_{\sigma(M)^{\prime}}(x)=tr_{\sigma(M)^{\prime}%
}(\Theta_{T}(x))+tr_{(\varphi(M)\otimes I)^{\prime}}(\Delta(I\otimes
x)\Delta).
\]
In particular, if $(T,\sigma)$ is an isometric representation of $\mathcal{E}$
then
\[
tr_{\sigma(M)^{\prime}}(\Theta_{T}(x))=dim_{l}(\mathcal{E})\cdot
tr_{\sigma(M)^{\prime}}(x)
\]
for every positive $x$ in $\sigma(M)^{\prime}$.
\end{corollary}

\begin{proof}
It follows from the computation in Proposition~\ref{estimate} that we have,
\[
tr_{\sigma(M)^{\prime}}(\Theta_{T}(x))+tr_{(\varphi(M)\otimes I)^{\prime}%
}((I\otimes c)\Delta^{2}(I\otimes c^{\ast}))=dim_{l}(\mathcal{E})\cdot
tr_{\sigma(M)^{\prime}}(x)
\]
where $\;c^{\ast}c=x\;$. Using the trace property of $tr_{(\varphi(M)\otimes
I)^{\prime}}$ we see that the second summand on the left hand side is equal to
$\;tr_{(\varphi(M)\otimes I)^{\prime}}(\Delta(I\otimes x)\Delta)\;$ and this
completes the proof of the main assertion. When the representation is
isometric, we have $\;\tilde{T}^{\ast}\tilde{T}=I\;$ and $\Delta=0\;$; thus
the conclusion follows.
\end{proof}

\begin{remark}
\label{remdim}Using the notation of Proposition~\ref{estimate} the result of
the corollary can also be written as
\[
dim_{l}(\mathcal{E})=\frac{tr_{\mathcal{B}^{\prime}}\left(
\begin{array}
[c]{cc}%
\Theta_{T}(x) & 0\\
0 & 0
\end{array}
\right)  +tr_{\mathcal{B}^{\prime}}\left(
\begin{array}
[c]{cc}%
0 & 0\\
0 & \Delta(I\otimes x)\Delta
\end{array}
\right)  }{tr_{\mathcal{B}^{\prime}}\left(
\begin{array}
[c]{cc}%
x & 0\\
0 & 0
\end{array}
\right)  }%
\]
for any positive element $x$ in $\sigma(M)^{\prime}$ whose trace is finite .
(We can, in fact, let $tr_{\mathcal{B}^{\prime}}$ be any (faithful, normal,
semifinite ) trace on $\mathcal{B}^{\prime}$).
\end{remark}

\subsection{Dilations of representations of correspondences}

We conclude our preliminary discussion by recalling and developing some
results concerning the minimal isometric dilation $(V,\rho)$ of a
representation $(T,\sigma)$ of a correspondence $\mathcal{E}$ over a von
Neumann algebra $M$. For details beyond those presented here see
\cite{MStensor} and \cite{QMP}.

Assume that $(T,\sigma)$ is a representation of $\mathcal{E}$ on $H$ and let
$\tilde{T}$ be the contraction defined above. Write $\Delta=(I-\tilde{T}%
^{\ast}\tilde{T})^{1/2}$ (in $B(\mathcal{E}\otimes_{\sigma}H)$) and let
$\mathcal{D}$ be the closure of the range of $\Delta$. Owing to formula
(\ref{Ttilde}), $\Delta$ commutes with $\sigma^{\mathcal{E}}\circ\varphi$ and
so $\mathcal{D}$ reduces $\sigma^{\mathcal{E}}\circ\varphi$. We let
$\sigma_{1}:=\sigma^{\mathcal{E}}\circ\varphi|\mathcal{D}$. Also let $L_{\xi
}:H\rightarrow\mathcal{E}\otimes_{\sigma}H$ be the map defined by $\;L_{\xi
}h=\xi\otimes h\;$ and write $\;D(\xi)=\Delta\circ L_{\xi}$. Note that
$\;T(\xi)=\tilde{T}\circ L_{\xi}$. The representation space $K$ of $(V,\rho)$
is
\[
K=H\oplus\mathcal{D}\oplus(\mathcal{E}\otimes_{\sigma_{1}}\mathcal{D}%
)\oplus(\mathcal{E}^{\otimes2}\otimes_{\sigma_{1}}\mathcal{D})\oplus\ldots
\]
Let $\sigma_{n+1}=\sigma_{1}^{\mathcal{E}^{\otimes n}}\circ\varphi_{n}$,
$n\geq1$, so that $\sigma_{n+1}$ is a representation of $M$ on $\mathcal{E}%
^{\otimes n}\otimes_{\sigma_{1}}\mathcal{D}$. The representation $\rho$ of $M$
on $K$ that we want is $\rho=\sigma\oplus\sum_{n\geq1}^{\oplus}\sigma_{n}$.
The map $V:\mathcal{E}\rightarrow B(K)$ is defined by
\[
V(\xi)=\left(
\begin{array}
[c]{cccc}%
T(\xi) & 0 & 0 & \ldots\\
D(\xi) & 0 & 0 & \ldots\\
0 & L_{\xi} & 0 & \\
0 & 0 & L_{\xi} & \\
&  &  & \ddots
\end{array}
\right)
\]
where $L_{\xi}$ here is the obvious map from $\mathcal{E}^{\otimes m}%
\otimes_{\sigma_{1}}\mathcal{D}$ to $\mathcal{E}^{\otimes(m+1)}\otimes
_{\sigma_{1}}\mathcal{D}$. Letting $\tilde{V}:\mathcal{E}\otimes_{\rho
}K\rightarrow K$ be the map sending $\xi\otimes k$ to $V(\xi)k$ , we can
write
\begin{equation}
\tilde{V}=\left(
\begin{array}
[c]{cccc}%
\tilde{T} & 0 & 0 & \ldots\\
\Delta & 0 & 0 & \ldots\\
0 & I & 0 & \\
0 & 0 & I & \\
&  &  & \ddots
\end{array}
\right)  \label{Vtilde}%
\end{equation}
where the identity operators in this matrix should be interpretted as the
operators that identify the spaces $\mathcal{E}\otimes_{\sigma_{n+1}%
}(\mathcal{E}^{\otimes n}\otimes_{\sigma_{1}}\mathcal{D})$ with $\mathcal{E}%
^{\otimes(n+1)}\otimes_{\sigma_{1}}\mathcal{D}$.

We shall refer to this $(V,\rho)$ as the \emph{minimal isometric dilation} of
$(T,\sigma)$ as it satisfies the following properties.

\begin{enumerate}
\item[(1)] $(V,\rho)$ is an isometric covariant representation of
$\mathcal{E}$ on $K$,

\item[(2)] $H$ reduces $\rho$ , and $\rho(a)|H=P_{H}\rho(a)|H=\sigma(a)$ for
all $a$ in $M$,

\item[(3)] $H^{\bot} = K \ominus H $ is invariant under each $V(\xi)$ ,
$\xi\in\mathcal{E}$ ; i.e., $P_{H}V(\xi)|H^{\bot} =0$,

\item[(4)] $P_{H}V(\xi)|H = T(\xi)$ , for all $\xi\in\mathcal{E}$, and

\item[(5)] the smallest subspace of $K$ containing $H$ and invariant under any
$V(\xi)$ is all of $K$.
\end{enumerate}

It is shown in \cite[Proposition 3.2]{MStensor} that a minimal isometric
dilation is unique up to unitary equivalence. Hence we refer to the
representation defined above as \emph{the} minimal isometric representation of
$(T,\sigma)$.

>From the isometric property of $(V,\rho)$ it follows that $\tilde{V}$ is an
isometry and, thus, the map $\Theta_{V}$, which, recall, is defined on
$\rho(M)^{\prime}$, is in fact a $\ast$-endomorphism of $\rho(M)^{\prime}$.
\ By Proposition and Definition \ref{theta}, $\Theta_{V}$ is unital precisely
when $V$ is fully coisometric and this happens if and only if $\tilde{V}$ is a
coisometry. \ However, a simple calculation shows that $\tilde{V}$ is a
coisometry if and only if $T$ is a coisometry, i.e., by Proposition and
Definition \ref{theta} again, if and only if $\Theta_{T}$ is unital.

We shall write $P_{n}$, $n\geq0$, for the projection $\Theta_{V}^{n}(I)$
($=\tilde{V}_{n}\tilde{V}_{n}^{\ast}\;$). Then the sequence $\{P_{n}%
\}_{n\geq0}$ is a decreasing sequence of projections and we write $P_{\infty
}:=\wedge P_{n}$. Moreover, if $Q_{n}:=P_{n}-P_{n+1}$,$\;n\geq0$, then the
$Q_{n}$ are mutually orthogonal (they are ``wandering'' projections) and
$\sum_{k=0}^{\infty}Q_{k}=I-P_{\infty}$.

The following terminology comes from \cite{wold}.

\begin{definition}
\label{induced}An isometric representation $(V,\rho)$ of $\mathcal{E}$ is said
to be \emph{induced} if there is a (normal) representation $\pi_{0}$ on a
Hilbert space $H_{0}$ such that $(V,\rho)$ is unitarily equivalent to the
representation $(R,\sigma)$ defined on $\mathcal{F}(\mathcal{E})\otimes
_{\pi_{0}}H_{0}$ by the formulae $\sigma=\pi_{0}^{\mathcal{F}(\mathcal{E}%
)}\circ\varphi_{\infty}$ and $R(\xi)=\pi_{0}^{\mathcal{F}(\mathcal{E})}%
(T_{\xi})$, $\xi\in\mathcal{E}$; i.e.,
\[
\sigma(a)=\varphi_{\infty}(a)\otimes I_{H_{0}},\;\;a\in M,
\]
and
\[
R(\xi)=T_{\xi}\otimes I_{H_{0}},\;\;\xi\in\mathcal{E},
\]
where, recall, $T_{\xi}$ is the creation operator on $\mathcal{F}%
(\mathcal{E})$ determined by $\xi$ and $\varphi_{\infty}$ is the diagonal
representation of $M$ on $\mathcal{F}(\mathcal{E})$ defined in the discussion
following Definition~\ref{correspondence}.
\end{definition}

We state for reference a result from \cite{wold} that will prove useful in our
analysis; it is a generalization of the Wold decomposition of an isometry.

\begin{proposition}
\label{decomp}\cite[Theorem 2.9]{wold} Every isometric representation
$(V,\rho)$ of $\mathcal{E}$ on $K$ decomposes as the direct sum
\begin{equation}
(V,\rho)=(V_{ind},\rho_{ind})\oplus(V_{\infty},\rho_{\infty})
\label{wolddecomp}%
\end{equation}
where $(V_{\infty},\rho_{\infty})$ is the restriction of $(V,\rho)$ to
$P_{\infty}(K)$ and is fully coisometric (and isometric), while $(V_{ind}%
,\rho_{ind})$ is the restriction of $(V,\rho)$ to $K\ominus P_{\infty}(K)$ and
is an induced representation. In fact, the representation $\pi_{0}$ of $M$
appearing in the definition of induced representation is, in this case, the
restriction of $\rho$ to the range of $Q_{0}$. (Here $P_{\infty}$ and $Q_{0}$
are as defined above).
\end{proposition}

\begin{definition}
\label{pure}The direct sum decomposition of an isometric representation of
$\mathcal{E}$ , $(V,\rho)$, given in equation (\ref{wolddecomp}) is called
\emph{the Wold decomposition} of $(V,\rho)$; $(V_{ind},\rho_{ind})$ is called
the \emph{pure part }or \emph{induced part} of $(V,\rho)$ and $(V_{\infty
},\rho_{\infty})$ is called the \emph{fully coisometric part} or
\emph{residual part }of $(V,\rho)$. If $(V_{\infty},\rho_{\infty})$ reduces to
zero, we say that $(V,\rho)$ is \emph{pure} or \emph{induced}.

If $(V,\rho)$ is the minimal isometric dilation of a representation
$(T,\sigma)$ of $\mathcal{E}$, we will say that $(T,\sigma)$ is \emph{pure }if
$(V,\rho)$ is pure.
\end{definition}

The ``purity'' of $T$ is reflected in $\Theta_{T}$ and $\Theta_{V}$ as the
next proposition indicates.

\begin{proposition}
\label{pureEquiv}Let $(T,\sigma)$ be a representation of $\mathcal{E}$ on a
Hilbert space $H$ and let $(V,\rho)$ be its minimal isometric dilation acting
on the Hilbert space $K$ containing $H$. Then the following conditions are equivalent.

\begin{enumerate}
\item[(1)] $(T,\sigma)$ is pure.

\item[(2)] $\Theta_{V}^{k}(I)\rightarrow0\;$ in the strong operator topology
on $K$.

\item[(3)] $\Theta_{T}^{k}(I)\rightarrow0\;$ in the strong operator topology
on $H$.
\end{enumerate}
\end{proposition}

\begin{proof}
The equivalence of conditions (1) and (2) is proved in Corollary 2.10 of
\cite{wold}. Condition (2) implies condition (3) simply because $\tilde{T}%
_{k}\tilde{T}_{k}^{\ast}=P_{H}\tilde{V}_{k}\tilde{V}_{k}^{\ast}P_{H}\;$. For
the other direction note that, since $\Vert\tilde{T}_{k}^{\ast}h\Vert
=\Vert\tilde{V}_{k}^{\ast}h\Vert$ for all $h$ in $H$, condition (3) implies
that $H$ is orthogonal to the range of $P_{\infty}=\wedge\tilde{V}_{k}%
\tilde{V}_{k}^{\ast}\;$. Hence it follows from the minimality of the dilation
that $P_{\infty}=0$.
\end{proof}

\section{The Curvature of a Representation}

In this section we define and study the curvature invariant for
representations of a given $W^{\ast}$-correspondence $\mathcal{E}$ over a
semifinite factor $M$.

\begin{definition}
\label{curvature}Let $M$ be a semifinite factor and let $\mathcal{E}$ be a
$W^{\ast}$-correspondence over $M$ with finite left dimension $d=dim_{l}%
(\mathcal{E})$. For every representation $(T,\sigma)$ of $\mathcal{E}$ we
define the \emph{curvature of }$(T,\sigma)$ to be
\begin{equation}
K(T,\sigma,\mathcal{E})=lim_{k\rightarrow\infty}\frac{tr_{\sigma(M)^{\prime}%
}(I-\Theta_{T}^{k}(I))}{\sum_{j=0}^{k-1}d^{j}} \label{curv}%
\end{equation}
where and $\Theta_{T}$ is the contractive, normal, completely positive map
defined in Proposition~\ref{theta}.
\end{definition}

When $\sigma$ and $\mathcal{E}$ are fixed, we sometimes abbreviate
$K(T,\sigma,\mathcal{E})$ as $K(T)$.

In order to prove the existence of the limit in the definition of the
curvature we require the following ``summability result'', which stated by
Popescu in \cite[p. 280]{gP01}. It is reminiscent of the fact that if a
sequence is convergent, then the sequence of its arithmetic means converges to
its limit.

\begin{lemma}
\label{limit}Let $\{a_{j}\}_{j=0}^{\infty}$ and $\{b_{j}\}_{j=0}^{\infty}$ be
two sequences of real numbers such that $b_{j}>0$ for all $j\geq0$. Write
$A_{k}$ and $B_{k}$ for the partial sums, $A_{k}=\sum_{j=0}^{k-1}a_{j}$ and
$B_{k}=\sum_{j=0}^{k-1}b_{j}$, and assume that $B_{k}\rightarrow\infty$ as
$k\rightarrow\infty$. If the limit $L=\lim_{j\rightarrow\infty}a_{j}/b_{j}$
exists (and finite) , then
\[
L=\lim_{k\rightarrow\infty}\frac{A_{k}}{B_{k}}.
\]
\end{lemma}

\begin{theorem}
\label{Curvature}Let $M$ be a semifinite factor, let $\mathcal{E}$ be a
$W^{\ast}$-correspondence over $M$, with finite left dimension $d=dim_{l}%
(\mathcal{E})$, and let $(T,\sigma)$ be a representation of $\mathcal{E}$ on a
Hilbert space. Then:

\begin{enumerate}
\item[(1)] The limit defining $K(T,\sigma,\mathcal{E})$ (in
Definition~\ref{curvature}) exists, either as a positive number or $+\infty$.

\item[(2)] $K(T,\sigma,\mathcal{E})=\infty\;$ if and only if $\;tr_{\sigma
(M)^{\prime}}(I-\Theta_{T}(I))=\infty$.

\item[(3)] If $tr_{\sigma(M)^{\prime}}(I-\Theta_{T}(I))\neq\infty$, so
$K(T,\sigma,\mathcal{E})<\infty$, then the following formulae hold.

\begin{enumerate}
\item[(3a)] If $d$ is larger than or equal to $1$, then
\[
K(T,\sigma,\mathcal{E})=lim_{k\rightarrow\infty}\frac{tr_{\sigma(M)^{\prime}%
}(\Theta_{T}^{k}(I)-\Theta_{T}^{k+1}(I))}{d^{k}}.
\]
and, if $d$ is strictly larger than $1$, then this limit is also given by the
equation $K(T,\sigma,\mathcal{E})=(d-1)lim_{k\rightarrow\infty}%
\frac{tr_{\sigma(M)^{\prime}}(I-\Theta_{T}^{k}(I))}{d^{k}}.$

\item[(3b)] If $d$ is strictly less than $1$, then $A:=\lim_{k\rightarrow
\infty}tr_{\sigma(M)^{\prime}}(I-\Theta_{T}^{k}(I))$ exists (and is finite)
and
\[
K(T,\sigma,\mathcal{E})=(1-d)A\text{.}%
\]
\end{enumerate}
\end{enumerate}
\end{theorem}

\begin{proof}
Write $tr$ for $tr_{\sigma(M)^{\prime}}$, $\Theta$ for $\Theta_{T}$ and
$a_{j}$ for $tr(\Theta^{j}(I)-\Theta^{j+1}(I))$. Then it follows from
Proposition~\ref{estimate} that
\begin{equation}
a_{j+1}\leq da_{j} \label{d}%
\end{equation}
for $j\geq0$. This shows that $a_{0}=\infty$ if and only if $a_{j}=\infty$ for
all $j\geq0$. Now assume $a_{0}<\infty$ and consider the sequence
$\{a_{j}/d^{j}\}_{j=0}^{\infty}$.It follows from equation (\ref{d}) that this
is a non-increasing sequence of nonnegative numbers . We write $L$ for its
limit. Then $0\leq L\leq a_{0}$.

Note that $tr(I-\Theta^{k}(I))=\sum_{j=0}^{k-1}a_{j}$.

Suppose $d\geq1$. Then we can use Lemma~\ref{limit} (with $b_{j}=d^{j}$) to
conclude that in this case, $K(T,\sigma,\mathcal{E})$ exists and equals $L$.
This proves one equality in part (3). For the other one note that
\[
\sum_{j=0}^{k-1}d^{j}=\frac{d^{k}-1}{d-1}%
\]
and $\lim_{k\rightarrow\infty}d^{k}/(d^{k}-1)=1$ when $d>1$. This proves the
assertions in (3a).

If $d<1$, the denominator in the definition of the curvature has a finite
limit $\;1/(1-d)\;$ as $k\rightarrow\infty$. Since it follows from equation
(\ref{d}) that $a_{j}\leq d^{j}a_{0}$ for all $j\geq0$, the numerator in the
definition of the curvature tends to a finite limit (denoted $A$ in (3b)).
This completes the proof of assertion (3b) and also the proofs of assertions
(1) and (2) because we just showed that, whenever $a_{0}\neq\infty$, the limit
defining the curvature exists and is finite.
\end{proof}

As we noted in the Introduction, this result ``captures'' Kribs's definition
of the curvature of a $d$-tuple \cite[Definition 2.4]{dK01} and Popescu's
Corollary 2.7 in \cite{gP01}.

\begin{lemma}
\label{isometric}Let $M$ be a finite factor, let $\mathcal{E}$ be a $W^{\ast}%
$-correspondence over $M$ and let $(T,\sigma)$ be a representation of
$\mathcal{E}$ on a Hilbert space. Then:

\begin{enumerate}
\item[(1)] The curvature of $T$, $K(T,\sigma,\mathcal{E})$, satisfies the
inequality
\[
K(T,\sigma,\mathcal{E})\leq tr_{\sigma(M)^{\prime}}(I-\tilde{T}\tilde{T}%
^{\ast})\text{,}%
\]
and if the representation is isometric, equality holds.

\item[(2)] Suppose $tr_{\sigma(M)^{\prime}}(I-\tilde{T}\tilde{T}^{\ast
})<\infty$. Then $K(T,\sigma,\mathcal{E})=tr_{\sigma(M)^{\prime}}(I-\tilde
{T}\tilde{T}^{\ast})$ if and only if, for every $j\geq1$,
\begin{equation}
{tr_{\sigma(M)^{\prime}}(\Theta_{T}^{j}(I)-\Theta_{T}^{j+1}(I))=dim_{l}%
(\mathcal{E})\cdot tr_{\sigma(M)^{\prime}}(\Theta_{T}^{j-1}(I)-\Theta_{T}%
^{j}(I)).} \label{0}%
\end{equation}
\end{enumerate}
\end{lemma}

\begin{proof}
Write $a_{j}=tr(\Theta^{j}(I)-\Theta^{j+1}(I))$ as in the proof of
Theorem~\ref{Curvature}. We have $\;a_{j+1}\leq da_{j}\;$ where $\;d=dim_{l}%
(\mathcal{E})\;$ and (using Corollary~\ref{isomest}) for an isometric
representation we have equality. Hence
\[
tr(I-\Theta^{k}(I))=\sum_{j=0}^{k-1}a_{j}\leq a_{0}\sum_{j=0}^{k-1}d^{j}%
\]
and equality holds for an isometric representation. This proves assertion (1).
In fact, it shows that equality holds whenever $a_{j+1}=da_{j}$ for all
$j\geq0$ and this proves one implication in assertion (2). For the other
implication, assume that $K(T,\sigma,\mathcal{E})=a_{0}$. If $d\geq1$ it
follows that $a_{0}$ is the limit of the decreasing sequence $\{a_{j}/d^{j}\}$
whose first term is $a_{0}$. Hence the sequence is constant. If $d<1$ then it
follows that $a_{0}=\sum a_{j}-\sum da_{j}=\sum(a_{j+1}-da_{j})+a_{0}$ and,
again, we get $a_{j+1}=da_{j}$ for all $j$.
\end{proof}

The notion of ``pure rank'' which we define next is motivated by the concept
first presented by Davidson, Kribs and Shpigel in \cite{DKS}. \ Their concept
coincides with ours in the case when the factor $M$ reduces to the scalars
$\mathbb{C}$.

\begin{definition}
\label{rank}Let $M$ be a semifinite factor, let $\mathcal{E}$ $\ $be a
$W^{\ast}$-correspondence over $M$, and let $(T,\sigma)$ be a representation
of $\mathcal{E}$ on a Hilbert space.

\begin{enumerate}
\item[(1)] For an element $a$ in $\sigma(M)^{\prime}$, we write $rank_{\sigma
(M)^{\prime}}(a)$ ( or, simply, $rank(a)$) for $tr_{\sigma(M)^{\prime}}(r(a))$
where $r(a)$ is the range projection of $a$.

\item[(2)] We define the \emph{pure rank} of $(T,\sigma)$ to be
\[
pure\;rank(T,\sigma)=rank_{\sigma(M)^{\prime}}(I-\tilde{T}\tilde{T}^{\ast}).
\]
\end{enumerate}
\end{definition}

The inequalities displayed in the next proposition generalize Theorem 2.7 of
\cite{dK01}.

\begin{proposition}
\label{purerank}Let $M$ be a semifinite factor and let $\mathcal{E}$ be
$W^{\ast}$-correspondence over $M$. If $(T,\sigma)$ is a representation of
$\mathcal{E}$ on a Hilbert space and if $(V,\rho)$ is its minimal isometric
dilation, then
\[
K(V,\rho,\mathcal{E})=tr_{\rho(M)^{\prime}}(I-\tilde{V}\tilde{V}^{\ast
})=pure\;rank(V,\rho)=pure\;rank(T,\sigma)
\]%
\[
\geq tr_{\sigma(M)^{\prime}}(I-\tilde{T}\tilde{T}^{\ast})\geq K(T,\sigma
,\mathcal{E}).
\]
\end{proposition}

\begin{proof}
Most of the asserted inequalities have already been proven. What needs to be
proved here is the equality of the pure ranks of $(T,\sigma)$ and $(V,\rho)$.
For this we shall analyze the construction of $(V,\rho)$. Let the Hilbert
space of $(T,\sigma)$ be $H$ and let the Hilbert space of $(V,\rho)$ be $K$
constructed in Section \ref{prelim}. Write the matrix of $\tilde{V}$ as in
equation (\ref{Vtilde}), write $P$ for the projection of $K$ onto $H$, and $Q$
for the projection $I-\tilde{V}\tilde{V}^{\ast}$. Note that $Q\in
\rho(M)^{\prime}$. Further, denote the span $\overline{span}\{V(\xi
)k\;:\;\xi\in\mathcal{E},\;\;k\in K\}$ by $L(K)$. Then, for $\xi$ in
$\mathcal{E}$ and $k$ in $K$, we have $QV(\xi)k=(I-\tilde{V}\tilde{V}^{\ast
})\tilde{V}(\xi\otimes k)=0$, since $\tilde{V}$ is an isometry. Hence $Q$
vanishes on $L(K)$. Since $(V,\rho)$ is minimal, we have $K=H\vee L(K)$ and,
thus, $Q(H)$ is dense in the range of $Q$. Setting $S=QP$ we conclude that the
range projection of $S$ is $Q$. We now turn to showing that the range
projection of $S^{\ast}=PQ$ is equal to the range projection of $I-\tilde
{T}\tilde{T}^{\ast}$. For this we use equation (\ref{Vtilde}) to compute
\begin{equation}
Q=I-\tilde{V}\tilde{V}^{\ast}=\left(
\begin{array}
[c]{cccc}%
I_{H}-\tilde{T}\tilde{T}^{\ast} & -\tilde{T}\Delta & 0 & \ldots\\
-\Delta\tilde{T}^{\ast} & I_{\mathcal{D}}-\Delta^{2} & 0 & \\
0 & 0 & 0 & \\
\vdots &  &  & \ddots
\end{array}
\right)  . \label{Q}%
\end{equation}
It follows from this that the range of $PQ$ is the space
\[
\{(I_{H}-\tilde{T}\tilde{T}^{\ast})h+\tilde{T}\Delta^{2}(\xi\otimes
f):\;f,h\in H,\;\;\xi\in\mathcal{E}\}.
\]
But $\tilde{T}\Delta^{2}=\tilde{T}(I-\tilde{T}^{\ast}\tilde{T})=(I-\tilde
{T}\tilde{T}^{\ast})\tilde{T}$. Hence the range of $PQ$ is contained in the
range of $I-\tilde{T}\tilde{T}^{\ast}$. Since the other containment is
obvious, the range projection of $PQ=S^{\ast}$ is the range projection of
$I-\tilde{T}\tilde{T}^{\ast}$. As the range projections of $S$ and $S^{\ast}$
are equivalent in $\rho(M)^{\prime}$, they have the same trace. More
precisely, this argument shows that they have the same trace with respect to
$tr_{\rho(M)^{\prime}}$. But then Lemma~\ref{trace} shows that the traces in
$\rho(M)^{\prime}$ and in $\sigma(M)^{\prime}$ coincide on the range
projection of $I-\tilde{T}\tilde{T}^{\ast}$ .
\end{proof}

We conclude this section with a generalization of \cite[Theorem 3.4]{dK01} and
\cite[Theorem 3.4]{gP01}.

\begin{theorem}
\label{critisom}Let $M$ be a semifinite factor, let $\mathcal{E}$ be a
correspondence over $M$ with finite left dimension and let $(T,\sigma)$ be a
pure representation of $\mathcal{E}$ with finite pure rank. Then $(T,\sigma)$
is isometric (hence, necessarily, an induced representation) if and only if
$K(T,\sigma,\mathcal{E})=pure\;rank(T,\sigma)$.
\end{theorem}

\begin{proof}
If the representation is isometric (and pure) it follows from
Proposition~\ref{pure} that it is induced. It also follows from
Proposition~\ref{purerank} that the equality holds.

So we now assume that the equality holds. It follows from
Proposition~\ref{purerank} that
\[
K(T,\sigma,\mathcal{E})=tr_{\sigma(M)^{\prime}}(I-\tilde{T}\tilde{T}^{\ast}).
\]
The argument in the proof of Theorem~\ref{Curvature} shows that this holds
only if, for every $j\geq1$,
\[
tr(\Theta^{j}(I-\Theta(I)))(=tr(\Theta^{j}(I)-\Theta^{j+1}(I)))=d^{j}\cdot
tr(I-\Theta(I))
\]
where $tr$ stands for $tr_{\sigma(M)^{\prime}}$, $\Theta$ is $\Theta_{T}$ and
$d$ is $dim_{l}(\mathcal{E})$. Note that the assumption we made that the pure
rank of the representation is finite means that the trace of $I-\Theta(I)$
and, thus, that the trace of $I-\Theta^{j}(I)$ is finite for all $j>0$. The
computation in the proof of Proposition~\ref{estimate}, applied to
$\;x_{j}=I-\Theta^{j}(I)\;$, now shows that, letting $c_{j}$ be $(I-\Theta
^{j}(I))^{1/2}$, the traces of $I\otimes x_{j}\;(=(I\otimes c_{j}^{2}))$ and
of $\;(I\otimes c_{j})\tilde{T}^{\ast}\tilde{T}(I\otimes c_{j})$ (both
elements of $(\varphi(M)\otimes I)^{\prime}$) are equal. Since $\tilde{T}$ is
a contraction and the two elements have the same finite trace, it follows that
the elements are equal. Thus
\[
(I\otimes c_{j})(I-\tilde{T}^{\ast}\tilde{T})(I\otimes c_{j})=0,\;j>0\text{.}%
\]
Since $\Theta^{j}(I)\rightarrow0$ in the strong operator topology (as
$(T,\sigma)$ was assumed to be pure), we conclude that $c_{j}\rightarrow0\;$in
the strong operator topology. Hence $\tilde{T}^{\ast}\tilde{T}=I$; i.e.,
$(T,\sigma)$ is isometric.
\end{proof}

\section{The Index and Curvature for Completely Positive Maps}

We shall now apply the curvature invariant for representations of
correspondences to the study of contractive, completely positive maps on
semifinite factors. We start by describing how one can associate, to a
contractive completely positive map $\Theta$ on a von Neumann algebra $N$, a
$W^{\ast}$-correspondence $\mathcal{E}$ and a representation $(T,\sigma)$ of
$\mathcal{E}$ such that $\;\Theta=\Theta_{T}$, where $\Theta_{T}$ is the
completely positive map associated to the representation $(T,\sigma)$ as in
Proposition and Definition~\ref{theta}.

The construction we describe was presented in \cite{QMP} and details can be
found there. We call attention to the fact that in \cite{QMP} we assumed that
the completely positive map $\Theta$ under discussion is unital. However, the
arguments and conclusions, with minor modifications, only require that
$\Theta$ be contractive (as well as normal and completely positive).

So, let $N$ be a von Neumann algebra acting on a Hilbert space $H$ and let
$\Theta:N\rightarrow N\;$ be a contractive, normal and completely positive
map. View $\Theta$ also as a completely positive map from $N$ into $B(H)$ and
write $N\otimes_{\Theta}H$ for the space of its Stinespring dilation (see
\cite{stinespring}). Recall that $N\otimes_{\Theta}H$ is obtained from the
algebraic tensor product $N\otimes H$ via the process of completion using the
sesquilinear form $\langle\cdot,\cdot\rangle$ defined by the formula
\[
\langle a_{1}\otimes h_{1},a_{2}\otimes h_{2}\rangle=\langle h_{1}%
,\Theta(a_{1}^{\ast}a_{2})h_{2}\rangle,
\]
where $a_{i}\otimes h_{i}\in N\otimes H$. The Hausdorff completion of
$N\otimes H$ is a Hilbert space which will be denoted $N\otimes_{\Theta}H$.
The Stinespring representation of $N$ on this space is given by the formula
\[
\pi(a)(b\otimes h)=ab\otimes h,
\]
for $a\in N$ and $b\otimes h\in N\otimes_{\Theta}H$. Also the formula
\[
W_{\Theta}(h)=I\otimes h\text{,}\;\;h\in H\text{,}%
\]
defines a bounded operator $W_{\Theta}$ mapping $H$ into $N\otimes_{\Theta}H$
satisfying the equations $\Vert W_{\Theta}\Vert=\Vert\Theta(I)\Vert$ and
\[
\Theta(a)=W_{\Theta}^{\ast}\pi(a)W_{\Theta}\;,\;\;a\in N.
\]
One also finds, in particular, that $\;W_{\Theta}^{\ast}(a\otimes
h)=\Theta(a)h$. The $W^{\ast}$-correspondence of importance to us is the
following one.

\begin{proposition and definition}
\label{esubtheta} Let $N$ be a von Neumann algebra acting on the Hilbert space
$H$ and let $\Theta$ be a contractive, normal, completely positive map on $N$.
We define $\mathcal{E}_{\Theta}$ to be the space of all operators from $H$ to
$N\otimes_{\Theta}H$ that intertwine the identity representation of $N$ on $H$
and the Stinespring representation $\pi$ just described; i.e.,
\[
\mathcal{E}_{\Theta}=\{X:H\rightarrow N\otimes_{\Theta}H\;:\;Xa=\pi
(a)X\;,\;a\in N\}.
\]
Then $\mathcal{E}_{\Theta}$ is a $W^{\ast}$-correspondence over $N^{\prime}$ (
the commutant of $N$ on $H$), where the right and the left actions are defined
by the formulae
\[
X\cdot b=X\circ b\;\;,\;\;\varphi(b)X=(I\otimes b)\circ X\;,\;\;b\in
N^{\prime}\text{,}%
\]
and the $N^{\prime}$-valued inner product is given by the formula
\[
\langle X,Y\rangle=X^{\ast}Y\;,\;\;X,Y\in\mathcal{E}_{\Theta}.
\]
We call $\mathcal{E}_{\Theta}$ the \emph{Arveson-Stinespring correspondence
}associated to $\Theta$.
\end{proposition and definition}

\begin{proof}
This is proved in \cite[Proposition 2.3]{QMP} under the assumption that
$\Theta$ is unital. \ However, a moment's reflection reveals that this
assumption plays no material role in the construction of $\mathcal{E}_{\Theta
}$.
\end{proof}

\begin{definition}
\label{idenrep}Let $\Theta$ be a contractive, normal, completely positive map
on the von Neumann algebra $N$ acting on the Hilbert space $H$ and let
$\mathcal{E}_{\Theta}$ be the Arveson-Stinespring correspondence associated to
$\Theta$. We define a representation $(T,\sigma)$ of $\mathcal{E}_{\Theta}$ on
$H$ by letting $\sigma$ be the identity representation of $N^{\prime}$ on $H$
and by setting
\begin{equation}
T(X)h=W_{\Theta}^{\ast}Xh,\;\;X\in\mathcal{E}_{\Theta},\;h\in H.
\label{identityrep}%
\end{equation}
It is shown in \cite[Discussion following Theorem and Definition 2.18]{QMP}
that this is indeed a representation of $\mathcal{E}_{\Theta}$, which will be
called the \emph{identity representation }of $\mathcal{E}_{\Theta}$.
\end{definition}

For future reference, we record the following fact that justifies the
terminology. It is essentially \cite[Corollary 2.23]{QMP}. The proof presented
there works here as well, and so will be omitted.

\begin{proposition}
\label{equality}Let $N$ be a von Neumann algebra acting on a Hilbert space $H$
and let $\Theta$ be a contractive, normal, completely positive map acting on
$N$. If $(T,\sigma)$ is the identity representation of the Arveson-Stinespring
correspondence $\mathcal{E}_{\Theta}$, then $\Theta=\Theta_{T}$.
\end{proposition}

Our next major objective, Theorem \ref{uniqued}, is to show that the left
dimension, $\dim_{l}(\mathcal{E}_{\Theta})$, of the Arveson-Stinespring
correspondence associated with a contractive, normal, completely positive map
$\Theta$ on a semifinite factor $N$ is independent of the Hilbert space on
which $N$ is represented. Along the way, we shall make some useful auxiliary observations.

First observe that if $N\otimes_{\Theta}H$ is the Hilbert space of the
Stinespring dilation of a contractive, normal, completely positive map
$\Theta$ acting on a von Neumann algebra $N\subseteq B(H)$, and if $a$ is an
operator in $N^{\prime}$, then the operator $I\otimes_{\Theta}a$ , acting on
$N\otimes_{\Theta}H$ according to the formula $\;(I\otimes_{\Theta}a)(b\otimes
h)=b\otimes ah\;$ is well defined and bounded, with $\;\Vert I\otimes_{\Theta
}a\Vert\leq\Vert a\Vert$. In fact, the map $a\rightarrow I\otimes_{\Theta}a$
is a normal representation of $N^{\prime}$ and we have the following result
that identifies this representation with the induced representation
$\sigma^{\mathcal{E}_{\Theta}}\circ\varphi_{\mathcal{E}_{\Theta}}$.

\begin{lemma}
\label{nthetah}If $\mathcal{E}_{\Theta}$ is the Arveson-Stinespring
correspondence associated with a contractive, normal, completely positive map
$\Theta$ on a von Neumann algebra $N$, then the operator $u:\mathcal{E}%
_{\Theta}\otimes H\rightarrow N\otimes_{\Theta}H\;$ defined by
\[
u(X\otimes h)=X(h)\;,\;X\in\mathcal{E}_{\Theta},\;h\in H
\]
is a Hilbert space isomorphism that intertwines the actions of $N^{\prime}$ on
the two spaces; i.e.,
\[
u(\sigma^{\mathcal{E}_{\Theta}}\circ\varphi_{\mathcal{E}_{\Theta}}(a))u^{\ast
}=u(\varphi_{\mathcal{E}_{\Theta}}(a)\otimes I)u^{\ast}=I\otimes_{\Theta
}a\;,\;\;a\in N^{\prime}.
\]
\end{lemma}

\begin{proof}
We compute:
\[
\langle X\otimes h,X\otimes h\rangle=\langle h,X^{\ast}Xh\rangle=\langle
X(h),X(h)\rangle.
\]
Hence $u$ is an isometry. The fact that $u$ is surjective follows from
\cite[Lemma 2.10]{QMP}. The intertwining property is a simple computation:
\[
u(\varphi(a)\otimes I)(X\otimes h)=u(\varphi(a)X\otimes h)=(\varphi
(a)X)(h)=X\otimes ah
\]%
\[
=(I\otimes a)(X\otimes h)=(I\otimes a)u(X\otimes h)
\]
for $X$ $\in\mathcal{E}_{\Theta}$ , $h\in H$ and $a\in N^{\prime}$, where we
have written $\varphi$ for $\varphi_{\mathcal{E}_{\Theta}}$.
\end{proof}

Let $\mathcal{C}=\{\left(
\begin{array}
[c]{cc}%
a & 0\\
0 & I\otimes_{\Theta}a
\end{array}
\right)  \in B(H\oplus(N\otimes_{\Theta}H)):\;a\in N^{\prime}\}$, a subalgebra
of $B(H\oplus(N\otimes_{\Theta}H))$ that is isomorphic to $N^{\prime}$. Let
$(T,\sigma)$ be the identity representation of $\mathcal{E}_{\Theta}$ (see
(\ref{identityrep})) and let $\mathcal{B}$ be the algebra defined in formula
(\ref{diag1}), associated with this representation (with $N^{\prime}$ playing
the role of $M$ there). Then the operator $\;U=I\oplus u$ implements an
isomorphism between $\mathcal{C}$ and $\mathcal{B}$. Explicitly, we have%
\begin{equation}
U\left(
\begin{array}
[c]{cc}%
a & 0\\
0 & \varphi(a)\otimes I_{H}%
\end{array}
\right)  U^{\ast}=\left(
\begin{array}
[c]{cc}%
a & 0\\
0 & I\otimes_{\Theta}a
\end{array}
\right)  \label{BandC}%
\end{equation}
for $a$ in $N^{\prime}$. (We suppressed the $\sigma$ in this formula.) Recall
that $\tilde{T}$ is the map from $\mathcal{E}_{\Theta}\otimes H$ to $H$
defined by $\tilde{T}(X\otimes h)=T(X)h=W_{\Theta}^{\ast}X(h)\;$. We shall
write $S$ for the map $S=\tilde{T}u^{\ast}\;$ and then we have,
\[
Su(X\otimes h)=\tilde{T}(X\otimes h)=W_{\Theta}^{\ast}X(h)
\]
and, if $X(h)$ may be expressed as $X(h)=b\otimes h\;$, then,
\[
S(b\otimes h)=W_{\Theta}^{\ast}(b\otimes h)=\Theta(b)h.
\]
Also,
\begin{equation}
U\left(
\begin{array}
[c]{cc}%
0 & \tilde{T}\\
0 & 0
\end{array}
\right)  U^{\ast}=\left(
\begin{array}
[c]{cc}%
0 & S\\
0 & 0
\end{array}
\right)  . \label{TandS}%
\end{equation}
In particular, the matrix in the right hand side of this equation lies in
$\mathcal{C}^{\prime}$.

Given an element $b$ in $N$, the operator $I_{\mathcal{E}}\otimes b$ lies in
$B(\mathcal{E}_{\Theta}\otimes H)$ and for $X\otimes h$ in $\mathcal{E}%
_{\Theta}\otimes H$ we have $u(I_{\mathcal{E}}\otimes b)(X\otimes
h)=u(X\otimes bh)=X(bh)=(b\otimes I_{H})X(h)\;$ (where the last equality
follows from the intertwining property of $X$). Thus, $\;u(I\otimes b)u^{\ast
}=b\otimes I\;$ and
\begin{equation}
U\left(
\begin{array}
[c]{cc}%
0 & 0\\
0 & I\otimes b
\end{array}
\right)  U^{\ast}=\left(
\begin{array}
[c]{cc}%
0 & 0\\
0 & b\otimes I
\end{array}
\right)  . \label{Itensorb}%
\end{equation}

The proof of Theorem \ref{uniqued}, which shows that the left dimension of the
Arveson-Stinespring correspondence associated with a completely positive map
on a semifinite factor is independent of the Hilbert space on which the factor
is represented, will involve a few steps. We start by fixing a (faithful,
normal) representation $\pi$ of the semifinite factor $N$ (on which the
normal, contractive, completely positive map $\Theta$ is defined) and we write
$\mathcal{E}$ for $\mathcal{E}_{\Theta}$. We also write $H_{\infty}$ for the
direct sum of infinitely many copies of $H$ and $\pi_{\infty}$ for infinite
ampliation of $\pi$ acting on $H_{\infty}$. (When it is clear in the
discussion to follow which representation we have in mind, we will suppress
$\pi$ or $\pi_{\infty}$.) We shall write $\mathcal{E}_{\infty}$, $(T_{\infty
},\sigma_{\infty})$, $\mathcal{C}_{\infty}$, $S_{\infty}$, etc. for the
objects defined using $\pi_{\infty}$ and $H_{\infty}$ in place of $\pi$ and
$H$. Note that $(\pi_{\infty}(N))^{\prime}$ is isomorphic to $M_{\infty}%
(\pi(N)^{\prime})$ (i.e., $\pi(N)^{\prime}\otimes B(l_{2})$) and
$\sigma_{\infty}$ is its identity representation on $H_{\infty}$. Our aim is
to show that $\;dim_{l}(\mathcal{E})=dim_{l}(\mathcal{E}_{\infty})\;$.

Let $v$ be the map from $N\otimes_{\Theta}H_{\infty}\;$to $(N\otimes_{\Theta
}H)_{\infty}$ defined by the formula
\[
v(b\otimes(h_{i}))=(b\otimes h_{i})\text{,}%
\]
$b\otimes(h_{i})\in N\otimes_{\Theta}H_{\infty}$. Then $v$ evidently is a
Hilbert space isomorphism. We write $V$ for the Hilbert space isomorphism
$\;V=I\oplus v$ mapping $H_{\infty}\oplus(N\otimes_{\Theta}H_{\infty})$ to
$H_{\infty}\oplus(N\otimes_{\Theta}H)_{\infty}$. Given an operator $R$ from
one Hilbert space to another, we write $R^{(\infty)}$ for its infinite ampliation.

The following lemma is easy to check and so we omit the proof.

\begin{lemma}
\label{cinfty}With $V$ as above we have,
\[
V\mathcal{C}_{\infty}V^{\ast}=M_{\infty}(\mathcal{C})\text{,}%
\]
and
\[
V\mathcal{C}_{\infty}^{\prime}V^{\ast}=M_{\infty}(\mathcal{C})^{\prime}%
\cong\mathcal{C}^{\prime}\otimes I_{\infty}\cong\mathcal{C}^{\prime}\;.
\]
\end{lemma}

We shall identify $M_{\infty}(\mathcal{C})^{\prime}\;$ with $\;\mathcal{C}%
^{\prime}\otimes I_{\infty}\;$, write $\Psi$ for the isomorphism of
$\mathcal{C}^{\prime}$ onto $\mathcal{C}^{\prime}\otimes I_{\infty}\;$ and
write $\Phi$ for the isomorphism of $\mathcal{C}^{\prime}$ onto $\;\mathcal{C}%
_{\infty}^{\prime}\;$ defined by
\[
\Phi(R)=V^{\ast}\Psi(R)V\;,\;\;R\in\mathcal{C}^{\prime}.
\]

\begin{lemma}
\label{phi}For $b\in N$, $S$ defined above and for$\;D=(I-S^{\ast}S)^{1/2}\;$
and $D_{\infty}$ $=(I-S_{\infty}^{\ast}S_{\infty})^{1/2}$, we have the
following equations:

\begin{enumerate}
\item[(1)] $\Phi\left(
\begin{array}
[c]{cc}%
\pi(b) & 0\\
0 & 0
\end{array}
\right)  =\left(
\begin{array}
[c]{cc}%
\pi_{\infty}(b) & 0\\
0 & 0
\end{array}
\right)  \;,\;\;b\in N$,

\item[(2)] $\Phi\left(
\begin{array}
[c]{cc}%
0 & 0\\
0 & b\otimes I_{H}%
\end{array}
\right)  =\left(
\begin{array}
[c]{cc}%
0 & 0\\
0 & b\otimes I_{H_{\infty}}%
\end{array}
\right)  \;,\;\;b\in N$ ,

\item[(3)] $\Phi\left(
\begin{array}
[c]{cc}%
0 & S\\
0 & 0
\end{array}
\right)  =\left(
\begin{array}
[c]{cc}%
0 & S_{\infty}\\
0 & 0
\end{array}
\right)  ,\;\;$ and

\item[(4)] $\Phi\left(
\begin{array}
[c]{cc}%
0 & 0\\
0 & D
\end{array}
\right)  =\left(
\begin{array}
[c]{cc}%
0 & 0\\
0 & D_{\infty}%
\end{array}
\right)  .$
\end{enumerate}
\end{lemma}

\begin{proof}
The proof of equations (1) and (2) is straightforward. For equation (3), write
$S^{(\infty)}$ for the diagonal map from $(N\otimes_{\Theta}H)_{\infty}$ to
$H_{\infty}$ induced by $S$ and compute, for $h=(h_{i})$ in $H^{\infty}$ and
$b$ in $N$,
\[
S^{(\infty)}v(b\otimes h)=S^{(\infty)}((b\otimes h_{i}))=(S(b\otimes
h_{i}))=(\pi(\Theta(b))h_{i})
\]%
\[
=\pi_{\infty}(\Theta(b))h=S_{\infty}(b\otimes h).
\]
This proves equation (3). Equation (4) follows immediately from equation (3).
\end{proof}

\begin{lemma}
\label{d}Let $N$ and $\Theta$ be as above. Fix a (faithful, normal)
representation $\pi$ of $N$ on $H$ with $M=\pi(N)^{\prime}$ finite. Let
$\mathcal{E}$ (=$\mathcal{E}_{\Theta}$), $\mathcal{C}$, $S$ and $D$
(=$(I-S^{\ast}S)^{1/2}$) be as above and let $tr$ be any ( normal, faithful,
semifinite) trace on $\mathcal{C}^{\prime}$. Then, for every positive $b$ in
$N$ with finite trace (with respect to a semifinite trace on $N$),
\begin{equation}
dim_{l}(\mathcal{E})=\frac{tr\left(
\begin{array}
[c]{cc}%
\Theta(b) & 0\\
0 & 0
\end{array}
\right)  +tr\left(
\begin{array}
[c]{cc}%
0 & 0\\
0 & D(b\otimes I_{H})D
\end{array}
\right)  }{tr\left(
\begin{array}
[c]{cc}%
b & 0\\
0 & 0
\end{array}
\right)  }. \label{eqd}%
\end{equation}
\end{lemma}

\begin{proof}
The lemma follows from Remark~\ref{remdim} using equations (\ref{BandC}),
(\ref{TandS}) and (\ref{Itensorb}) (to translate from $\mathcal{B}$ to
$\mathcal{C}$).
\end{proof}

\begin{theorem}
\label{uniqued}Let $\Theta$ be a normal, contractive, completely positive map
on the semifinite factor $N$. Let $\pi_{1}$ and $\pi_{2}$ be two faithful,
normal $\ast$-representations of $N$ on Hilbert spaces $H_{1}$ and $H_{2}$
respectively such that both $\pi_{1}(N)^{\prime}$ and $\pi_{2}(N)^{\prime}$
are finite. Write $\mathcal{E}_{1}$ and $\mathcal{E}_{2}$ for the $W^{\ast}%
$-correspondences associated to $\Theta$ and the representations $\pi_{1}$ and
$\pi_{2}$ respectively. Then
\[
dim_{l}(\mathcal{E}_{1})=dim_{l}(\mathcal{E}_{2}).
\]
\end{theorem}

\begin{proof}
Using Lemma~\ref{d}, it suffices to show that the quotient on the right hand
side of equation (\ref{eqd}) is the same for $\pi_{1}$ and for $\pi_{2}$. We
shall show this for every $\pi_{1}$, $\pi_{2}$ (regardless of whether the
commutants of the images of these representations are finite or not). So we
now drop the assumption on the commutants.

If the two representations are equivalent then the conclusion of the theorem
certainly holds. This is the case if $\pi_{1}(N)^{\prime}$ and $\pi
_{2}(N)^{\prime}$ are both infinite (\cite[Exercises 9.6.30 and 9.6.31]{KR}).

Thus, it suffices to consider only the case where $\pi_{2}$ is $\pi_{1,\infty
}$. So we now write $\pi$ for $\pi_{1}$ and $\pi_{\infty}$ for $\pi_{2}$. But
then the the analysis above (Lemma~\ref{cinfty} and Lemma~\ref{phi}) is
applicable and the equality of the expressions on the right hand side of
equation (\ref{eqd}) follows from Lemma~\ref{phi}.
\end{proof}

\begin{definition}
\label{index}Let $N$ be a semifinite factor and let $\Theta$ be contractive,
normal, completely positive map on $N$. Then we define the \emph{index }of
$\ \Theta$ to be the left dimension of the Arveson-Stinespring correspondence
associated with any (necessarily faithful) normal representation of $N$ having
finite commutant. \ We denote the index of $\Theta$ by $d(\Theta)$.
\end{definition}

Theorem \ref{uniqued} shows that our notion of ``index'' is well defined. \ In
\cite{arvindex} Arveson defined a notion of index for continuous semigroups of
normal, contractive completely positive maps on $B(H)$. \ Our definition here
seems to be the analogous one for a single map on an arbitrary semifinite factor.

The following lemma gives a convenient computation of the index.

\begin{lemma}
\label{compind}Let $\Theta$ be a contractive, normal completely positive map
on the semifinite factor $N$. Fix a semifinite, normal, faithful trace $tr$ on
$N$ and a projection $e$ in $N$ with $tr(e)=1$ (if $N$ is finite we can assume
that $tr$ is normalized and $e=I$). Write $L^{2}(N)$ for $L^{2}(N,tr)$ and
denote by $\lambda$ and $\rho$ the left and right actions of $N$ on $L^{2}%
(N)$. Let $H$ be $\rho(e)L^{2}(N)$ (to be written $L^{2}(Ne)$) and let $M$ be
the algebra $\rho(eNe)$. Consider the Hilbert space $eN\otimes_{\Theta}H$ (a
subspace of $N\otimes_{\Theta}H$) on which $M$ acts according to the equation
\begin{equation}
\rho(eae)\cdot(eb\otimes h)=eb\otimes\rho(eae)h,\;\;a,b\in N,\;h\in H.
\label{action}%
\end{equation}
Then
\[
d(\Theta)=dim_{M}(eN\otimes H).
\]
\end{lemma}

\begin{proof}
Let $\pi$ be the restriction of the representation $\lambda$ of $N$ on
$L^{2}(N)$ to $H$. Then $M=\pi(N)^{\prime}\;$ and it is a finite factor. The
operator $u$ of Lemma~\ref{nthetah} maps $\mathcal{E}_{\Theta}\otimes H$ onto
$N\otimes_{\Theta}H$. Since $L^{2}(M)$ can be identified with $\lambda
(e)L^{2}(Ne)$, the space $\mathcal{E}_{\Theta}\otimes L^{2}(M)$ is
$(I_{\mathcal{E}}\otimes\lambda(e))(\mathcal{E}_{\Theta}\otimes H)$. As
$u(I_{\mathcal{E}}\otimes\lambda(e))u^{\ast}=e\otimes I\;$ (see equation
(\ref{Itensorb})), $\;u$ is a unitary operator mapping $\;\mathcal{E}_{\Theta
}\otimes L^{2}(M)\;$ onto $\;eN\otimes_{\Theta}H$. It follows from
Lemma~\ref{nthetah} that it is a module isomorphisms (when the latter space is
viewed as a module over $M$ as in equation(\ref{action})). Thus they have the
same dimension over $M$.
\end{proof}

\begin{proposition}
\label{ti}Let $N$ be a semifinite factor, let $t_{1},t_{2},\ldots t_{n}\;$ be
elements of $N$ such that $\sum t_{i}t_{i}^{\ast}\leq I$,$\;$ and $\Theta$ be
the contractive, normal, completely positive map on $N$ defined by the
formula
\[
\Theta(a)=\sum_{i=1}^{n}t_{i}at_{i}^{\ast}\text{.}%
\]
Then
\[
d(\Theta)=\dim_{\mathbb{C}}span\{t_{i}:\;1\leq i\leq n\}\in\{0,1,\ldots n\}
\]
where $\dim_{\mathbb{C}}$ is the dimension as a complex vector space.
\end{proposition}

\begin{proof}
By Lemma~\ref{compind} we have (once we fix a projection $e$ with trace equal
to 1),
\[
d(\Theta)=\dim_{\rho(eNe)}(eN\otimes_{\Theta}L^{2}(Ne)).
\]
Define the Hilbert space $K$ by the formula%
\[
K=\overline{span}\{\sum_{i}\oplus eat_{i}^{\ast}be\;:\;a,b\in N\;\}\subseteq
\sum_{i}\oplus L^{2}(eNe).
\]
In fact, $\sum\oplus L^{2}(eNe)\;$ is a module over $M=\rho(eNe)$ and $K$ is a
submodule. Define a linear operator $w:eN\otimes_{\Theta}L^{2}(Ne)\rightarrow
K\;$ by the formula $\;w(ea\otimes be)=\sum\oplus eat_{i}^{\ast}be$. This
defines $w$ on a dense subspace and we compute, for $a,b,c$ and $f$ in $N$,
\[
\langle\sum\oplus eat_{i}^{\ast}be,\sum\oplus ect_{i}^{\ast}fe\rangle
=\sum\langle eat_{i}^{\ast}be,ect_{i}^{\ast}fe\rangle=\sum\langle
be,t_{i}a^{\ast}ect_{i}^{\ast}fe\rangle
\]%
\[
=\langle be,\Theta(a^{\ast}ec)fe\rangle=\langle ea\otimes_{\Theta}%
be,ec\otimes_{\Theta}fe\rangle.
\]
Hence $w$ can be extended to a unitary operator (also denoted $w$). It is also
a module isomorphism and, thus, we have $d(\Theta)=dim_{M}(K)$. Now write $Q$
for the projection, in $B(\sum\oplus L^{2}(eNe))$, whose range is $K$.
Identifying $B(\sum\oplus L^{2}(eNe))$ with $B(L^{2}(eNe))\otimes
M_{n}(\mathbb{C})$, we see that $Q$ lies in the commutant of both
$\lambda(eNe)\otimes I$ and $\rho(eNe)\otimes I$. Hence $Q=I\otimes q\;$ for
some projection $q$ in $M_{n}(\mathbb{C})$. Clearly, $\;d(\Theta)=\dim(q)\;$
and, in particular, $d(\Theta)$ is an integer between $0$ and $n$. Write $k$
for $d(\Theta)=\dim(q)\;$ and let $m=\dim span\{t_{i}\}=\dim span\{t_{i}%
^{\ast}\}$. We can find partial isometries in $N$, $\{w_{l}\}$, with
$w_{l}w_{l}^{\ast}=e\;$ for all $l$ and $\sum w_{j}^{\ast}w_{j}=I$. (If $N$ is
finite , this set is just $\{I\}$.) For every choice of $l$ and $p$, the
vector $\sum_{i}\oplus w_{l}t_{i}^{\ast}w_{p}^{\ast}\;$ lies in $K$. We shall
write this vector as a column vector $(w_{l}t_{i}^{\ast}w_{p}^{\ast})$ and
then we have $(I-q)(w_{l}t_{i}^{\ast}w_{p})=0$. Since $\;t_{i}^{\ast}%
=\sum_{l,p}w_{l}^{\ast}(w_{l}t_{i}^{\ast}w_{p}^{\ast})w_{p}$, we find that the
(column) vector $(t_{i}^{\ast})$ satisfies $(I-q)(t_{i}^{\ast})=0$. As
$rank(I-q)=n-k$, we conclude that $m\;(=\dim span\{t_{i}^{\ast}\})\leq k$. For
the other inequality note that since $\dim span\{t_{i}^{\ast}\}=m$, we can
find an $(n-m)\times n$ matrix $A$ (over $\mathbb{C}$) with $rankA=n-m$ such
that $A(t_{i}^{\ast})=0$. But then $K$ lies in the kernel of $A\otimes I$ and,
consequently, $Aq=0$. We conclude that $m=n-(n-m)=\dim KerA\geq dim\;q=k$,
and, therefore, that $m=k$.
\end{proof}

In the next proposition, $D$ is the operator on $N\otimes_{\Theta}H$ described
in Lemma~\ref{d}.

\begin{proposition}
\label{scaleatx}Let $\Theta$ be a completely positive map with finite index
acting on a semifinite factor $N$ and let $x$ be a positive element of $N$
with finite trace. Then the following assertions are equivalent.

\begin{enumerate}
\item[(1)] $tr(\Theta(x))=d(\Theta) \cdot tr(x)$.

\item[(2)] $D(x\otimes I_{H})D = 0 $.

\item[(3)] The element $c=x^{1/2}$ satisfies the equation,
\[
\Theta(ac)\Theta(cb)=\Theta(ac^{2}b)\text{,}%
\]
for all $a$ and $b$ in $N$.
\end{enumerate}

Moreover, if any of these conditions holds for $x$ then they hold for all
positive $y$ in $xNx$ (in particular, they hold for all $y\in N$ such that
$0\leq y\leq x$).
\end{proposition}

\begin{proof}
The equivalence of parts (1) and (2) follows from Lemma~\ref{d}.

Suppose (2) holds and write $c=x^{1/2}$. Then $(c\otimes I)D=0$ and, thus,
$(c\otimes I)D^{2}(c\otimes I)=0$. Since $D^{2}=I-S^{\ast}S$, we have
$(c\otimes I)S^{\ast}S(c\otimes I)=c^{2}\otimes I$ and
\[
S(ac\otimes I)S^{\ast}S(cb\otimes I)S^{\ast}=S(ac^{2}b\otimes I)S^{\ast}%
\]
for all $a$ and $b$ in $N$. Recall that $S$ is the operator from
$N\otimes_{\Theta}H$ to $H$ (for a fixed representation of $N$ on $H$) defined
by $S(b\otimes h)=\Theta(b)h\;,\;b\in N\;,\;h\in H$. It is easy to check that
$S^{\ast}h=I\otimes h$ for $h$ in $H$. Thus, for $y$ in $N$,
\[
S(y\otimes I)S^{\ast}h=S(y\otimes I)(I\otimes h)=\Theta(y)h\;,\;\;h\in H.
\]
Thus, $S(y\otimes I)S^{\ast}=\Theta(y)$ and this completes the proof of (3).
For the other direction, condition (3) implies that $S(a^{\ast}c\otimes
I)D^{2}(ca\otimes I)S^{\ast}=0$ for all $a\in N$. Thus $D(c\otimes I)(a\otimes
I)S^{\ast}S=0$ for all $a\in N$. Note that the range of $S^{\ast}S$ is
$I\otimes\lbrack\Theta(N)(H)]$ and, in order to prove (2), it suffices to show
that the subspace of $N\otimes_{\Theta}H$ spanned by elements of the form
$(a\otimes I)(I\otimes\Theta(b)k)$ for $a$ and $b$ in $N$ and $k$ in $H$ is
dense. Suppose $\sum n_{i}\otimes h_{i}$ in $N\otimes_{\Theta}H$ is orthogonal
to this subspace. Then
\[
0=\sum\langle n_{i}\otimes h_{i},a\otimes\Theta(b)k\rangle=\sum\langle
\Theta(a^{\ast}n_{i})h_{i},\Theta(b)k\rangle.
\]
Since this holds for all $b\in N$ and $k\in H$,
\[
\sum\Theta(a^{\ast}n_{i})h_{i}=0.
\]
>From this we conclude that $\sum n_{i}\otimes h_{i}$ is orthogonal to
$a\otimes h$ for all $a\in N$ and $h\in H$; i.e., it is equal to $0$. This
shows that the subspace above is indeed dense and completes the proof of (2).

The last statement of the proposition is clear by considering condition (2).
\end{proof}

We now can define our principle object of study, the curvature of a completely
positive map.

\begin{definition}
\label{curvcp}Let $N$ be a semifinite factor with faithful, normal semifinite
trace, $tr$, and let $\Theta$ be a contractive, normal completely positive map
on $N$ with finite index $d=d(\Theta)$. We define the \emph{curvature} of
$(\Theta,tr)$, $K(\Theta,tr)$, by the formula
\[
K(\Theta,tr)=\lim_{k\rightarrow\infty}\frac{tr(I-\Theta^{k}(I))}{\sum
_{j=0}^{k-1}d^{j}}\text{.}%
\]
\end{definition}

Theorem~\ref{uniqued} and Theorem~\ref{Curvature} guaranty that the curvature
is well defined. Of course, the statements of Theorem~\ref{Curvature} have
obvious analogues for $K(\Theta,tr)$.

Observe that if $\Theta$ is a contractive, normal, completely positive map on
a von Neumann algebra, then the sequence $\{\Theta^{k}(I)\}_{k\geq0}$ is a
decreasing sequence of positive operators and so converges to a positive
operator in the strong operator topology. The following terminology is
consistent with Definition \ref{pure} and will be useful in the sequel.

\begin{definition}
\label{purefull}Let $\Theta$ be a contractive, normal, completely positive map
on a von Neumann algebra. We write $\Theta^{\infty}(I)$ for the limit of
\ $\{\Theta^{k}(I)\}_{k\geq0}$ in the strong operator topology. Then $\Theta$
will be called \emph{full} if $\Theta^{\infty}(I)=I$ and \emph{pure} if
$\Theta^{\infty}(I)=0$.
\end{definition}

\begin{theorem}
\label{k=a0}Let $\Theta$ be a normal, contractive, completely positive map
with finite index acting on a semifinite factor $N$. Then

\begin{enumerate}
\item[(1)] $K(\Theta,tr) \leq tr(I-\Theta(I)) .$

\item[(2)] $K(\Theta,tr)=\infty\;$ if and only if $\;tr(I-\Theta(I))=\infty$.

\item[(3)] Suppose $tr(I-\Theta(I))<\infty$. Then $K(\Theta,tr)=tr(I-\Theta
(I))$ if and only if, for all $a$ and $b$ in $N$,
\[
\Theta(ac_{\infty})\Theta(c_{\infty}b)=\Theta(ac_{\infty}^{2}b)
\]
where $c_{\infty}=(I-\Theta^{\infty}(I))^{1/2}$. In particular, if $\Theta$ is
pure, then the equality holds in the first assertion, (1), if and only if
$\Theta$ is a $\ast$-endomorphism.
\end{enumerate}
\end{theorem}

\begin{proof}
Assertion (1) follows from Lemma~\ref{isometric} and the second assertion
follows from assertion (2) of Theorem~\ref{Curvature}. Now assume that
$tr(I-\Theta(I))\;$ is finite. For every $k$ write $x_{k}=\Theta
^{k-1}(I)-\Theta^{k}(I)$. Then $c_{\infty}^{2}=\sum_{k=1}^{\infty}x_{k}$
(where the convergence is in the strong operator topology). The condition
\[
\Theta(ac_{\infty})\Theta(c_{\infty}b)=\Theta(ac_{\infty}^{2}b)
\]
for all $a$ and $b$ in $N$, is equivalent to condition (2) of
Proposition~\ref{scaleatx} for $x=c_{\infty}^{2}$. But this holds if and only
if it holds for all $x_{k}$. Using that proposition, this is equivalent to the
condition
\[
tr(\Theta(x_{k}))=d(\Theta)\cdot tr(x_{k})\;,\;\;k\geq1.
\]
Applying part (2) of Lemma~\ref{isometric}, we find that the latter condition
is equivalent to $K(\Theta,tr)=tr(I-\Theta(I))$.
\end{proof}

Theorem \ref{k=a0} suggests that, for non unital maps, it is better to
``normalize'' $K$ and study the normalized curvature of a completely positive map.

\begin{definition}
\label{k1}Let $\Theta$ be a normal, contractive completely positive non unital
map acting on a semifinite factor $N$ and assume that $\Theta$ has finite
index and finite curvature. Then the \emph{normalized curvature} of $\Theta$,
denoted $K_{1}(\Theta)$, is to be
\[
K_{1}(\Theta)=K(\Theta,tr)/tr(I-\Theta(I))
\]
where $tr$ is any semifinite, normal, faithful trace on $N$.
\end{definition}

With the definition of $K_{1}$ in hand, the following corollary of Theorem
\ref{k=a0} is immediate.

\begin{corollary}
\label{endom}Let $\Theta$ be a normal $\ast$-endomorphism of $N$ with finite
index. Then

\begin{enumerate}
\item[(1)] For all positive $x$ in $N$, $\;tr(\Theta(x))=d(\Theta)\cdot tr(x)$.

\item[(2)] $K(\Theta,tr)=tr(I-\Theta(I))$; hence either $\Theta(I)=I\;$ or
$\;K_{1}(\Theta)=1$.

\item[(3)] If $\Theta_{1}$ is another normal $\ast$-endomorphism of $N$ with
finite index, then $d(\Theta\circ\Theta_{1})=d(\Theta)d(\Theta_{1})$.
\end{enumerate}
\end{corollary}

Specializing still further, we have the following corollary, in which
$\operatorname{mod}(\alpha)$ denotes the modulus of an automorphism $\alpha$
of a type $II_{1}$ factor (See \cite{KR}.)

\begin{corollary}
\label{autom}Let $N$ be a factor of type $II_{1}$ with normalized trace $\tau
$, and let $\alpha$ be a normal $\ast$-automorphism of $N$. Then
$d(\alpha)=mod(\alpha)\;$and$\;K(\alpha,\tau)=0.$
\end{corollary}

\begin{example}
\label{p}Let $N$ be a type $II_{1}$ factor and $p$ be a projection such that
$pNp$ is isomorphic to $N$. Write $\Theta:N\rightarrow pNp\;$ for this
isomorphism and view $\Theta$ as an endomorphism on $N$. Then $d(\Theta
)=\tau(p)$ and $K(\Theta,\tau)=1-\tau(p)$, where $\tau$ is the normalized
trace on $N$.
\end{example}

This follows easily from Corollary~\ref{endom}, since $\Theta(I)=p$.

\vspace{5 mm}

\begin{example}
\label{exp}Let $N$ be a type $II_{1}$ factor with normalized trace $\tau$ and
let $N_{0}\subseteq N$ be a subfactor of finite index, $[N:N_{0}]$. If $E$ is
the (trace preserving) conditional expectation from $N$ onto $N_{0}$, then
$d(E)=[N:N_{0}]$ and$\;K(E,\tau)=0.$
\end{example}

\begin{proof}
The statement about $K$ is obvious since the map is unital. As for its index,
recall that, using Lemma~\ref{compind}, we have
\[
d(E)=dim_{\rho(N)}(N\otimes_{E}L^{2}(N)).
\]
Let $e_{0}=e_{N_{0}}$ be the projection of $L^{2}(N,\tau)$ onto $L^{2}%
(N_{0},\tau)\;$ (extending the map $E$). Let $N_{1}$ be the von Neumann
algebra generated by $N$ and $e_{0}$. It is the next algebra in the basic
construction of Jones. Write $\tau_{1}$ for the normalized trace of $N_{1}$.
Then $\;\tau_{1}(xe_{0})=[N:N_{0}]^{-1}\tau(x)\;$ for all $x\in N\;$
(\cite[Proposition 3.1.2]{JS}). We define a linear map $\;S:N\otimes_{E}%
L^{2}(N,\tau)\rightarrow L^{2}(N_{1},\tau_{1})\;$ by $\;S(a\otimes
_{E}b)=ae_{0}b\;$ for $a,b$ in $N$, and compute: For $a,b,c$ and $d$ in $N$,
\[
\langle a\otimes_{E}b,c\otimes_{E}d\rangle=\langle b,E(a^{\ast}c)d\rangle
=\tau(b^{\ast}E(a^{\ast}c)d)=\tau(db^{\ast}E(a^{\ast}c))
\]%
\[
=\tau_{1}(db^{\ast}E(a^{\ast}c)e_{0})[N:N_{0}]=\tau_{1}(db^{\ast}e_{0}a^{\ast
}ce_{0})[N:N_{0}]
\]%
\[
=\tau_{1}(b^{\ast}e_{0}a^{\ast}ce_{0}d^{\ast})[N:N_{0}]=[N:N_{0}]\langle
ae_{0}b,ce_{0}d\rangle_{L^{2}(N_{1},\tau_{1})}.
\]
Hence the map $\;V=[N:N_{0}]^{-1/2}S\;$ extends to a Hilbert space isomorphism
from $N\otimes_{E}L^{2}(N,\tau)$ onto $L^{2}(N_{1},\tau_{1})$. (The
surjectivity of $V$ uses \cite[2.6(d)]{J}). In fact, $V$ intertwines the
actions of $\rho(N)$:
\[
V(\rho(c)(a\otimes_{E}b))=V(a\otimes_{E}bc)=[N:N_{0}]^{-1/2}ae_{0}%
bc=\rho(c)V(a\otimes_{E}b)
\]
for $a,b$ and $c$ in $N$. Thus, if we write $N_{2}$ for the next algebra in
Jones's basic construction (containing $N_{1}$), we have
\[
d(E)=dim_{\rho(N)}L^{2}(N_{1},\tau_{1})=dim_{N_{2}^{\prime}}L^{2}(N_{1}%
,\tau_{1})=[N_{1}^{\prime}:N_{2}^{\prime}]dim_{N_{1}^{\prime}}L^{2}(N_{1}%
,\tau_{1})
\]%
\[
=[N_{1}^{\prime}:N_{2}^{\prime}]=[N_{2}:N_{1}]=[N:N_{0}].
\]
In the computation above we used the fact that, on $L^{2}(N_{1},\tau_{1})$,
$N_{2}^{\prime}=\rho(N)\;$ (see Proposition 3.1.2 (ii) of \cite{JS}) and the
equality $\;dim_{N_{2}^{\prime}}H=[N_{1}^{\prime}:N_{2}^{\prime}%
]dim_{N_{1}^{\prime}}H$, which holds for every $N_{1}^{\prime}$-module $H$.
\end{proof}

\vspace{5 mm}

In order to show that $d(\Theta)$ is an outer conjugacy invariant we need
first to recall some results from \cite{QMP}. Let $\Theta_{1}$ and $\Theta
_{2}$ be two normal, contractive, completely positive maps on $N$ and let $H$
be a left $N$-module. Recall the definition of $\mathcal{E}_{\Theta_{i}}$:
\[
\mathcal{E}_{\Theta_{i}}=\{X:H\rightarrow N\otimes_{\Theta_{i}}%
H\;:\;Xa=(a\otimes I)X\;,\;\;a\in N\}
\]
for $i=1,2$. We can also define the space $N\otimes_{\Theta_{1}}%
N\otimes_{\Theta_{2}}\otimes H\;$ as the Hausdorff completion of the space we
get when we define the inner product
\[
\langle a\otimes_{\Theta_{1}}b\otimes_{\Theta_{2}}h,c\otimes_{\Theta_{1}%
}d\otimes_{\Theta_{2}}k\rangle=\langle h,\Theta_{2}(b^{\ast}\Theta_{1}%
(a^{\ast}c)d)k\rangle
\]
on the algebraic tensor product. We now set
\[
\mathcal{E}=\{X:H\rightarrow N\otimes_{\Theta_{1}}N\otimes_{\Theta_{2}%
}H\;:\;Xa=(a\otimes I\otimes I)X\;,\;\;a\in N\}.
\]
Define the map $\;\Psi:\mathcal{E}_{\Theta_{2}}\otimes\mathcal{E}_{\Theta_{1}%
}\rightarrow\mathcal{E}\;$ by
\[
\Psi(X\otimes Y)=(I\otimes X)Y
\]
where $I\otimes X$ is the map from $N\otimes_{\Theta_{1}}H$ to $N\otimes
_{\Theta_{1}}N\otimes_{\Theta_{2}}H$ given by the equation $(I\otimes
X)(a\otimes h)=a\otimes Xh$. We also define a map $V_{0}:N\otimes_{\Theta
_{2}\Theta_{1}}H\rightarrow N\otimes_{\Theta_{1}}N\otimes_{\Theta_{2}}H$ via
the equation $V_{0}(a\otimes h)=a\otimes I\otimes h$ and we let $V:\mathcal{E}%
_{\Theta_{2}\Theta_{1}}\rightarrow\mathcal{E}$ be the map defined by
$\;V(X)=V_{0}X$.

The following proposition can be found in \cite[Propositions 2.12 and
2.14]{QMP}. It was presented there under the assumption that the maps
$\Theta_{1}$ and $\Theta_{2}$ are unital, but the proof given holds without
this assumption.

\begin{proposition}
\label{mpsi}In the notation just established, we have the following.

\begin{enumerate}
\item[(1)] The map $\Psi$ is an isomorphism of correspondences.

\item[(2)] The map $V$ is an isometric bimodule map.

\item[(3)] The map $m=V^{\ast}\Psi$ is a coisometry and $m^{\ast}%
:\mathcal{E}_{\Theta_{2}\Theta_{1}}\rightarrow\mathcal{E}_{\Theta_{2}}%
\otimes\mathcal{E}_{\Theta_{1}}\;$ is an isometric bimodule map.

\item[(4)] If either $\Theta_{2}$ is an endomorphism (of $N$) or $\Theta_{1}$
is an automorphism, then the map $m$ is an isomorphism of correspondences.
Hence, in either case,
\[
\mathcal{E}_{\Theta_{2}\Theta_{1}}\simeq\mathcal{E}_{\Theta_{2}}%
\otimes\mathcal{E}_{\Theta_{1}}.
\]
\end{enumerate}
\end{proposition}

Combining this proposition with Corollary~\ref{multdim} we immediately get the
following inequality. (Compare with Corollary~\ref{endom}(3) ).

\begin{proposition}
\label{submult}For two normal, contractive, completely positive maps
$\Theta_{1}$ and $\Theta_{2}$ on a semifinite factor $N$, we have
\[
d(\Theta_{1}\Theta_{2})\leq d(\Theta_{1})d(\Theta_{2}).
\]
\end{proposition}

We also obtain

\begin{theorem}
\label{invariant}Let $\Theta_{1}$ and $\Theta_{2}$ be two normal, contractive,
completely positive maps on a semifinite factor $N$.

\begin{enumerate}
\item[(1)] If the maps are outer conjugate (i.e., if there is an automorphism
$\alpha$ of $N$ and a unitary operator $u$ in $N$ such that $\;ad(u)\circ
\Theta_{1}=\alpha^{-1}\circ\Theta_{2}\circ\alpha$), then $d(\Theta
_{1})=d(\Theta_{2})$. In particular, one has a finite index if and only if the
other one has.

\item[(2)] If they both have a finite index and they are conjugate (i.e., if
$\Theta_{1}=\alpha^{-1}\circ\Theta_{2}\circ\alpha$ for some automorphism
$\alpha$ on $N$) then $K(\Theta_{1},tr)=K(\Theta_{2},tr)d(\alpha)$, where $tr$
is any faithful normal trace on $N$ and where, recall, $d(\alpha)$ is the
index of $\alpha$.$\;$Further, if the curvature is finite and the maps are non
unital, then $K_{1}(\Theta_{1})=K_{1}(\Theta_{2})$.
\end{enumerate}
\end{theorem}

\begin{proof}
If the maps are outer conjugate and we write $\beta$ for $ad(u)$, it follows
from Proposition~\ref{mpsi} that
\[
\mathcal{E}_{\alpha}\otimes\mathcal{E}_{\beta}\otimes\mathcal{E}_{\Theta_{1}%
}\simeq\mathcal{E}_{\Theta_{2}}\otimes\mathcal{E}_{\alpha}.
\]
Using Corollary~\ref{multdim} we find that
\[
d(\alpha)d(\beta)d(\Theta_{1})=d(\Theta_{2})d(\alpha)
\]
and, since $d(\beta)=1$ (by Example~\ref{ti}), the assertions in part (1)
follow. The assertions in part (2) now result from the following computation.
\[
tr(I-\Theta_{2}^{k}(I))=tr(I-\alpha\circ\Theta_{1}^{k}(\alpha^{-1}%
(I)))=tr(\alpha(I-\Theta_{1}^{k}(I)))
\]%
\[
=d(\alpha)tr(I-\Theta_{1}^{k}(I)).
\]
\end{proof}

Given a completely positive map $\Theta$ on $N$, we showed in \cite{QMP} that
we can construct a ``dilation'' of $\Theta$ to a $\ast$-endomorphism $\alpha$
on another von Neumann algebra $R$ in which $N$ ``sits'' as a corner $R$. We
turn to describing this construction and relating the invariants of $\Theta$
to those of $\alpha$. So, let $\Theta$ be a contractive, normal, completely
positive map on the semifinite factor $N$ and let $N$ act on the Hilbert space
$H$. Form the Arveson-Stinespring correspondence $\mathcal{E}_{\Theta}$ and
its identity representation $(T,\sigma)$ and recall that $\Theta=\Theta_{T}$
by Proposition \ref{equality}. Let $(V,\rho)$ be the minimal isometric
dilation of $(T,\sigma)$ and write $K$ for its representation space (see the
discussion preceding Definition~\ref{induced}). We let $R$ be $\rho(N^{\prime
})^{\prime}$ and we let $\alpha$ be the endomorphism of $R$, $\Theta_{V}$. If
we write $W$ for the isometric embedding of $H$ onto $K$, then $N$ is
$W^{\ast}RW$ and $\alpha$ is a ``dilation'' of $\Theta$ in the sense portrayed
in the following proposition which can be found in \cite[Theorem 2.24]{QMP}.

\begin{proposition}
\label{dilation}With the notation just established, we may assert that:

\begin{enumerate}
\item[(1)] $W^{\ast}RW=N$ and $R^{\prime}$ is a normal homomorphic image of
$N^{\prime}$.

\item[(2)] $\alpha$ is a normal $\ast$-endomorphism of $R$.

\item[(3)] For every non-negative integer $n$ ,
\[
\Theta^{n}(a)=W^{\ast}\alpha^{n}(WaW^{\ast})W
\]
and
\[
\Theta^{n}(W^{\ast}bW)=W^{\ast}\alpha(b)W
\]
for all $a\in N$ and $b\in R$.
\end{enumerate}
\end{proposition}

It should be noted that this proposition was proved in \cite{QMP} under the
assumption that $\Theta$ is unital and the proof presented there seems to use
this fact by virtue of \cite[Theorem 2.20]{QMP}, which requires that $\Theta$
be unital. However, one can use instead \cite[Theorem and Definition
2.18]{QMP} directly, which does not require that $\Theta$ be unital. \ The
dilation endomorphism $\alpha$ that is produced, when $\Theta$ is not unital,
is not unital either.

The dilation $\alpha$ of $\Theta$ is unique up to conjugacy (\cite[Theorem
5.4]{wAp01}) and so really is an artifact of $\Theta$, independent of how $N$
is represented.

\begin{theorem}
\label{invofdil}Let $\Theta$ be a contractive, normal, completely positive map
on the semifinite factor $N$ and let $\alpha$ be the dilation of $\Theta$
constructed above from the identity representation $(T,\sigma)$ of the
Arveson-Stinespring correspondence $\mathcal{E}_{\Theta}$ associated with
$\Theta$.

\begin{enumerate}
\item[(1)] If $\Theta$ has finite index, then so does $\alpha$ and
$d(\Theta)=d(\alpha)$.

\item[(2)] If $d(\Theta)<\infty$, then $K(\alpha,tr_{\rho(N^{\prime})^{\prime
}})\geq K(\Theta,tr_{\sigma(N^{\prime})^{\prime}})$.

\item[(3)] If $d(\Theta)<\infty$ and $\Theta$ is pure, then $K(\alpha
,tr_{\rho(N^{\prime})^{\prime}})=tr_{\rho(N^{\prime})^{\prime}}(I-\alpha(I)).$

\item[(4)] If $d(\Theta)<\infty$, if $\Theta$ is pure and if $K(\alpha
,tr_{\rho(N^{\prime})^{\prime}})=K(\Theta,tr_{\sigma(N^{\prime})^{\prime}})$,
then $\Theta=\alpha$.
\end{enumerate}
\end{theorem}

\begin{proof}
Write $\mathcal{E}$ for $\mathcal{E}_{\Theta}$ and $\mathcal{E}_{V}$ for
$\mathcal{E}_{\Theta_{V}}=\mathcal{E}_{\alpha}$ (where $(V,\rho)$ is the
isometric representation of $\mathcal{E}$ dilating $(T,\sigma)$). In order to
prove that $d(\alpha)=d(\Theta)$ note that $\mathcal{E}$ is a $W^{\ast}%
$-correspondence over $N^{\prime}$ while $\mathcal{E}_{V}$ is a $W^{\ast}%
$-correspondence over $\rho(N^{\prime})$. We shall construct a map
$\;U:\mathcal{E}\rightarrow\mathcal{E}_{V}\;$ satisfying

\begin{enumerate}
\item[(i)] $U$ is a linear bijection,

\item[(ii)] it is a bimodule map in the sense that, for every $X$ in
$\mathcal{E}$ and $a,b$ in $N^{\prime}$,
\[
U(aXb)=\rho(a)U(X)\rho(b)
\]
and,

\item[(iii)] $\langle U(X_{1}),U(X_{2}) \rangle= \rho(\langle X_{1},X_{2}
\rangle)$.
\end{enumerate}

Once we construct such a map, assertion (1) will follow. To define $U$ we let
$u:R\otimes_{\alpha}K\rightarrow K\;$ be the linear map defined by the
equation $u(a\otimes k)=\alpha(a)k\;$. It is easy to check that $u$ is an
isometry onto $[\alpha(R)K]$. Also note that, for every $b\otimes g$ in
$R\otimes_{\alpha}K$ and every $k$ in $K$, we have $\;\langle u^{\ast
}k,b\otimes g\rangle=\langle k,\alpha(b)g\rangle=\langle I\otimes k,b\otimes
g\rangle\;$; hence $\;u^{\ast}k=I\otimes k$. It also follows from this that
$\;u^{\ast}\rho(a)=(I\otimes\rho(a))u^{\ast}\;$ for $a$ in $N^{\prime}$. Now
set
\[
U(X)=u^{\ast}V(X)\;\;,\;\;X\in\mathcal{E}.
\]
For $X$ in $\mathcal{E}$, $U(X)$ is a map from $K$ to $R\otimes_{\alpha}K$.
For $a$ in $R$ and $k$ in $K$, we have
\[
U(X)ak=u^{\ast}V(X)ak=u^{\ast}\tilde{V}(X\otimes ak)=u^{\ast}\tilde
{V}(I\otimes a)\tilde{V}^{\ast}\tilde{V}(X\otimes k)
\]%
\[
=u^{\ast}\alpha(a)V(X)k=a\otimes V(X)k=(a\otimes I)(I\otimes V(X)k)
\]%
\[
=(a\otimes I)u^{\ast}V(X)k=(a\otimes I)U(X)k.
\]
Hence $U$ maps $\mathcal{E}$ into $\mathcal{E}_{V}$.

For $a,b$ in $N^{\prime}$ and $X$ in $\mathcal{E}$ we have
\[
U(aXb)=u^{\ast}V(aXb)=u^{\ast}\rho(a)V(X)\rho(b)
\]%
\[
=(I\otimes\rho(a))u^{\ast}V(X)\rho(b)=\rho(a)\cdot U(X)\rho(b).
\]
Thus (ii) holds. To prove (iii) we compute (using the fact that $uu^{\ast
}=\alpha(I)=\tilde{V}\tilde{V}^{\ast}$)
\[
\langle U(X_{1}),U(X_{2})\rangle=U(X_{1})^{\ast}U(X_{2})=V(X_{1})^{\ast
}uu^{\ast}V(X_{2})
\]%
\[
=V(X_{1})^{\ast}V(X_{2})=\rho(\langle X_{1},X_{2}\rangle.
\]
This proves (iii). It is now clear that $U$ is injective and it is left to
prove surjectivity. For this assume that $Y$ in $\mathcal{E}_{V}$ is
orthogonal to the range of $U$. Then $Y^{\ast}u^{\ast}V(X)k=0$ for all $k$ in
$K$ and $X$ in $\mathcal{E}$. But the vectors of the form $V(X)k=\tilde
{V}(X\otimes k)\;$span a dense subspace of the range of $\tilde{V}$ which is
equal to the range of $u$. Since $u$ is an isometry, vectors of the form
$u^{\ast}V(X)k$ span a dense subspace of $K$. Thus $Y=0$ and, since
$\mathcal{E}_{V}$ is self dual, this implies that $U$ is surjective completing
the proof of (i),(ii) and (iii). Thus $\;d(\Theta)=dim_{l}(\mathcal{E}%
)=dim_{l}(\mathcal{E}_{V})=d(\alpha)$. Assertion (2) now follows from
Proposition~\ref{purerank}. Assertions (3) and (4) are immediate from Theorem
\ref{critisom}. Alternatively, use Theorem \ref{k=a0}.
\end{proof}

In \cite{Parrott} Parrott studied the curvature of a single contraction $T$.
His work was supplemented by Levy in \cite{rLpp00}. In our terminology the
curvature of a single contraction is the curvature of the contractive, normal,
completely positive map $\Theta(S)=TST^{\ast},$ for $S$ in $B(H)$. Our
analysis allows us to replace $B(H)$ by a general semifinite factor and to
extend Parrott's results.

\begin{proposition}
\label{singlet}Let $N$ be a semifinite factor with a normal semifinite trace
$tr$. Let $t\in N$ be a non zero contraction with $tr(I-tt^{\ast})<\infty$ and
set $\Theta(a)=tat^{\ast}\;$ for $a$ in $N$. Then

\begin{enumerate}
\item[(1)] $d(\Theta)=1$,

\item[(2)]
\[
K(\Theta,tr)=lim_{k\rightarrow\infty}\frac{tr(I-t^{k}t^{k\ast})}%
{k}=lim_{k\rightarrow\infty}tr(t^{k}t^{k\ast}-t^{k+1}t^{(k+1)\ast})
\]%
\[
=tr(I-tt^{\ast})-tr((I-t^{\ast}t)^{1/2}(I-t_{\infty})(I-t^{\ast}t)^{1/2})
\]
where $t_{\infty}=\lim t^{k}t^{k\ast}$ in the strong operator topology.

\item[(3)] If $\Theta$ is pure (i.e., if $t_{\infty}=0$) then
\[
K(\Theta,tr)=tr(I-tt^{\ast})-tr(I-t^{\ast}t).
\]
\end{enumerate}
\end{proposition}

\begin{proof}
The fact that $d(\Theta)=1\;$ follows from Example~\ref{ti}. Let $H$ be an
$N$-module as in Lemma~\ref{d}. Then the map $v:N\otimes_{\Theta}H\rightarrow
H\;$ defined by $\;v(a\otimes_{\Theta}h)=at^{\ast}h\;$ is a unitary operator.
The map $S$ of Lemma~\ref{phi} is defined by the equation $\;S(a\otimes
h)=\Theta(a)h=tat^{\ast}h$. Hence $Sv^{\ast}$ in $B(H)$ satisfies $Sv^{\ast
}(at^{\ast}h)=tat^{\ast}h$. Thus $Sv^{\ast}=t\;$ and, consequently,
$vDv^{\ast}=(I-t^{\ast}t)^{1/2}\;$. Also, for $b$ in $N$, $\;v(b\otimes
I_{H})v^{\ast}=b$. It now follows from Lemma~\ref{d} that, for a positive $b$
in $N$ with finite trace,
\[
tr(b)=d(\Theta)tr(b)=tr(\Theta(b))+tr((I-t^{\ast}t)^{1/2}b(I-t^{\ast}%
t)^{1/2}).
\]
Thus, for $x=I-\Theta(I)$ and $j\geq0$, we have
\[
tr(\Theta^{j}(x))-tr(\Theta^{j+1}(x))=tr((I-t^{\ast}t)^{1/2}\Theta
^{j}(x)(I-t^{\ast}t)^{1/2}).
\]
Summing this up for $j$ from $0$ to $k$ we find that
\[
tr(x)=tr(\Theta^{k+1}(x))+tr((I-t^{\ast}t)^{1/2}(x+\ldots+\Theta
^{k}(x))(I-t^{\ast}t)^{1/2}).
\]
But $x+\ldots+\Theta^{k}(x)=I-\Theta^{k+1}(I)$. Hence, when we take the limit,
as $k\rightarrow\infty$, and use the fact that $tr(\Theta^{k+1}(x))\rightarrow
K(\Theta,tr)$ (Theorem~\ref{Curvature}(3)), we get part (2). Part (3) follows
from (2).
\end{proof}

\end{document}